\documentclass[a4paper,11pt]{article}
\usepackage[utf8]{inputenc}  
\usepackage[OT2,T1]{fontenc}
\usepackage{babel}

\usepackage{amsthm}
\usepackage{amsmath}
\usepackage{amssymb}
\usepackage{mathrsfs}
\usepackage{dsfont}

\usepackage{enumitem}
\usepackage{ifthen}
\usepackage{verbatim}
\usepackage{geometry}
\usepackage{hyperref}
\usepackage{url}
\usepackage{fancyhdr}
\usepackage{tikz-cd}
\usepackage{stmaryrd}
\usepackage{authblk}
\usepackage{mathtools}
\usepackage{zref}

\newcommand{\R}{{\mathbb R}}

\newcommand{\Z}{{\mathbb Z}}
\newcommand{\Q}{{\mathbb Q}}
\renewcommand{\C}{{\mathbb C}}
\newcommand{\F}{{\mathbb F}}
\DeclareSymbolFont{cyrletters}{OT2}{wncyr}{m}{n}
\DeclareMathSymbol{\Sha}{\mathalpha}{cyrletters}{"58}
\newcommand{\rar}{\rightarrow}

\newcommand{\mrm}[1]{\mathrm{#1}}
\newcommand{\mfk}[1]{\mathfrak{#1}}
\newcommand{\OO}{\mathcal{O}}
\newcommand{\Qbar}{\overline{\mathbb{Q}}}
\newcommand{\Fr}{\mrm{Frob}}
\newcommand{\Aut}[1]{\ensuremath{\operatorname{Aut}(#1)}}
\newcommand{\Auti}[2]{\ensuremath{\operatorname{Aut}_{#1}(#2)}}
\newcommand{\GL}[1]{\ensuremath{\operatorname{GL}_2(#1)}}
\newcommand{\SL}[1]{\ensuremath{\operatorname{SL}_2(#1)}}
\newcommand{\PGL}[1]{\ensuremath{\operatorname{PGL}_2(#1)}}

\DeclareMathOperator{\Tr}{Tr}

\makeatletter
\let\oldlabel\label
\renewcommand{\label}[1]{%
  \zref@labelbylist{#1}{special}
  \oldlabel{#1}
}

\zref@newlist{special}
\zref@newprop{section}{\arabic{section}.\arabic{subsection}}
\zref@addprop{special}{section}

\newcounter{propo}
\renewcommand{\thepropo}{\thesection.\arabic{propo}}
\makeatletter \@addtoreset{propo}{section}\makeatother

\makeatletter \@addtoreset{rema}{section}\makeatother

\newcounter{theop}
\renewcommand{\thetheop}{\Alph{theop}}

\newcounter{qnc}

\newcommand{\lem}[2][None]{\refstepcounter{propo}\ifthenelse{\equal{#1}{None}}{}{\label{#1}} \smallskip \noindent \textbf{Lemma \thepropo\;} {\it #2}\smallskip}
\newcommand{\prop}[2][None]{\refstepcounter{propo}\ifthenelse{\equal{#1}{None}}{}{\label{#1}} \smallskip \noindent \textbf{Proposition \thepropo\;} {\it #2} \smallskip}
\newcommand{\propn}[2][None]{\refstepcounter{theop}\ifthenelse{\equal{#1}{None}}{}{\label{#1}}\smallskip \noindent \textbf{Theorem \thetheop\;} {\it #2} \smallskip}
\newcommand{\qnn}[2][None]{\refstepcounter{qnc}\ifthenelse{\equal{#1}{None}}{}{\label{#1}}\smallskip \noindent \textbf{Question \theqnc\;} {\it #2} \smallskip}
\newcommand{\cor}[2][None]{\refstepcounter{propo}\ifthenelse{\equal{#1}{None}}{}{\label{#1}} \smallskip \noindent \textbf{Corollary \thepropo\;} {\it #2} \smallskip}
\newcommand{\defi}[1]{\smallskip \noindent {\bf Definition\;} {#1} \smallskip}
\newcommand{\rem}[1]{\smallskip \noindent {\bf Remark\;} {#1} \smallskip}
\newcommand{\rems}[1]{\smallskip \noindent {\bf Remarks\;} {#1} \smallskip}
\newcommand{\demo}[1]{\noindent {\it Proof.\;--\;} #1\hfill$\Box$ \smallskip}

\binoppenalty=\maxdimen
\relpenalty=\maxdimen

\title{On obstructions to the Euler system method for Rankin-Selberg convolutions}
\author[]{Elie Studnia}

\begin{document}

\maketitle

\section*{Introduction}
Let $k,N \geq 1$ be integers. We denote by $\mathcal{S}_k(\Gamma_1(N))$ the space of cuspidal holomorphic modular forms with level $\Gamma_1(N)$. Let $f\in \mathcal{S}_k(\Gamma_1(N))$ be a cuspidal newform with character $\varepsilon_f$, and $L$ be a number field containing its coefficients. Given a maximal ideal $\mathfrak{p} \subset \OO_L$, we can attach to $f$ a free $\OO_{L_{\mathfrak{p}}}$-module $T_{f,\mfk{p}}$ of rank two endowed with a continuous linear action of $G_{\Q}:=\mrm{Gal}(\Qbar/\Q)$. By works of Momose \cite{Momose} and Ribet \cite{Ribet2}, it is known that, when $f$ has no complex multiplication and $k \geq 2$, the image of $G_{\Q}$ in $\Aut{T_{f,\mfk{p}}} \simeq \GL{\OO_{L_{\mathfrak{p}}}}$ contains a conjugate of $\SL{\Z_p}$ for all but finitely many $\mathfrak{p}$. More recently, Loeffler \cite{bigimage} proved an \emph{adelic} open image theorem for such newforms analogous to that of Serre \cite{Serre-image-ouverte} for elliptic curves without complex multiplication. 

Let $g \in \mathcal{S}_l(\Gamma_1(M))$ be another newform with character $\varepsilon_g$; after enlarging $L$ if necessary, we assume that its Fourier coefficients are contained in $L$. As above, we attach to $g$ a free $\OO_{L_{\mathfrak{p}}}$-module of rank two $T_{g,\mfk{p}}$ endowed with a continuous linear action of $G_{\Q}$. Let $T_{f,g,\mfk{p}} := T_{f,\mfk{p}} \otimes_{\OO_{L_{\mathfrak{p}}}} T_{g,\mfk{p}}$. The other goal of \cite{bigimage} was to study the image of $G_{\Q}$ in $\Aut{T_{f,g,\mfk{p}}}$. In particular, Loeffler proves that, when the following conditions are satisfied:

\begin{itemize}[noitemsep,label=\tiny$\bullet$]
\item $k,l \geq 2$ and $k$ or $l$ is odd\footnote{Loeffler's results apply somewhat more generally, but this is a simple sufficient condition.},
\item neither $f$ nor $g$ has complex multiplication,
\item $f$ and $g$ are not twists of each other,
\item $\varepsilon_f\varepsilon_g \neq 1$,
\end{itemize}
 then, for all but finitely many $\mathfrak{p}$, there exists a $\sigma \in G_{\Q(\mu_{p^{\infty}})} := \mrm{Gal}(\Qbar/\Q(\mu_{p^{\infty}}))$ such that $T_{f,g,\mfk{p}}/(\sigma-1)T_{f,g,\mfk{p}}$ is free of rank one. Loeffler also discusses the case where $k \geq 2, l=1$ and asks the following question: 

\qnn[loefflerqn]{(Loeffler, see \cite[Remark 4.4.2]{bigimage}) Assume that $k \geq 2$, $l=1$, $f$ does not have complex multiplication, and $\varepsilon_f\varepsilon_g \neq 1$. Is it true that, for all but finitely many $\mathfrak{p}$, there exists $\sigma \in G_{\Q(\mu_{p^{\infty}})}$ such that $T_{f,g,\mfk{p}}/(\sigma-1)T_{f,g,\mfk{p}}$ is free of rank one?}

Question \ref{loefflerqn} arises in the context of \emph{Euler systems}, a theory that aims to bound Selmer groups of $\mfk{p}$-adic Galois representations using certain classes in the Galois cohomology of its dual: a typical result is \cite[Theorem 5.2.2]{MR-ES}. For now, it seems that applying this method to a $\mfk{p}$-adic representation $T$ of $G_{\Q}$ requires the existence of $\sigma \in G_{\Q(\mu_{p^{\infty}})}$ such that $T/(\sigma-1)T$ is free of rank one over the coefficient ring. For instance, in \emph{loc.cit.}, this is assumption (H.2). For the Galois representations considered in this article, the Euler system would typically be the Euler system of Beilinson--Flach elements constructed in \cite{KLZ15}. 

Loeffler proves in \cite[Proposition 4.4.1]{bigimage} that the answer to Question \ref{loefflerqn} is positive when $N$ and $M$ are coprime. We prove in this article that the answer to Question \ref{loefflerqn} is positive if certain weaker conditions are satisfied. It should be stressed that these assumptions are not verified in general, but they are reasonably easy to check in practice and apply quite broadly. 

Before stating them, we need to recall the notion of inner twist, whose properties are stated more expansively in Section \ref{sect-setup}. 

Let $f \in \mathcal{S}_k(\Gamma_1(N))$ be a newform with $q$-expansion $f = \sum_{n\geq 1}{a_n(f)q^n}$: it is well-known that there is a number field $L \subset \C$ containing $a_n(f)$ for every $n \geq 1$. Let $\chi$ be a primitive Dirichlet character. We say that $f$ has an \emph{inner twist by $\chi$} if there is a field embedding $\sigma: L \rar \C$ and a positive multiple $N'$ of $N$ such that, for any integer $n \geq 1$ coprime to $N'$, $a_n(f)\chi(n) = \sigma(a_n(f))$. It can be shown (see Lemma \ref{momose15}) that, in such a case, the conductor of $\chi$ divides $N$, we can take $N'=N$, and one has $\sigma(L) = L$. 

\propn[weakerconds]{Let $\mathcal{C}$ denote the collection of primitive Dirichlet characters $\chi$ such that $f$ has an inner twist by $\chi$. The answer to Question \ref{loefflerqn} is positive when any of the following conditions holds:
\begin{itemize}[noitemsep,label=$-$]
\item For any $\chi \in \mathcal{C}$, $M$ and the conductor of $\chi$ are coprime. This condition is in particular satisfied when $f$ has no inner twist or when $N$ and $M$ are coprime. 
\item Every $\chi \in \mathcal{C}$ is even, i.e. satisfies $\chi(-1)=1$.
\item The group $\varepsilon_g\left(\bigcap_{\chi \in \mathcal{C}}{\ker{\chi}}\right)$ contains an element of order $2$ and $g$ does not have either real or complex multiplication. 
\item Let $\mathcal{K}$ denote the subgroup of $\hat{\Z}^{\times}$ given by $\mathcal{K}=\bigcap_{\chi \in \mathcal{C}}{\ker{\chi}}$. Then the group $\varepsilon_g\left(\mathcal{K}\right)$ contains an element of order $2$ but no element of order $4$, and, for every non-trivial Dirichlet character $\varepsilon$ such that $g \otimes \varepsilon=g$, the group $\varepsilon\left(\mathcal{K}\right)$ is not equal to $\{1\}$. 
\item $\varepsilon_f^2=\varepsilon_g^2=1$, and $\mathcal{C}$ does not contain $\varepsilon_g$ or any non-trivial quadratic Dirichlet character $\varepsilon$ such that $g \otimes \varepsilon = g$. 
\end{itemize}}

As explained above, because we can give a positive answer to Question \ref{loefflerqn} in more situations, we can apply the machinery of Euler systems with fewer technical assumptions. As a consequence, extending \cite[Theorem 11.7.4]{KLZ15}, we prove the following instance of the Bloch--Kato conjecture.

Recall the following definition first:

\defi{Let $E$ be an elliptic curve over a number field $K$, and $p$ be a prime number. The $p^{\infty}$-Selmer group of $E$ is the group 
\[\mrm{Sel}_{p^{\infty}}(E/K) = \ker\left[H^1(K,E[p^{\infty}]) \rar \prod_v{H^1(K_v,E(\overline{K_v}))}\right],\] where $v$ runs over the finite places of $K$ and the map is induced by localizing at the place $v$ and the inclusion $E[p^{\infty}] \subset E$. It is an abelian group where every element is killed by a power of $p$, so it has the natural structure of a $\Z_p$-module.  }

\propn[twistedbsd]{Let $E$ be an elliptic curve over $\Q$ without complex multiplication and $\rho$ be an odd two-dimensional irreducible representation of $G_{\Q}$, with splitting field $K$ and coefficients in some number field $L_0$ with ring of integers $\OO$. Let $\mathfrak{p}$ be a maximal ideal of $\OO$ with residue characteristic $p$. Suppose that the following assumptions hold:
\begin{enumerate}[noitemsep]
\item $p$ is coprime to $30N_{\rho}N_E$, where $N_E$ (resp. $N_{\rho}$) is the conductor of $E$ (resp. $\rho$), 
\item if $\rho$ has dihedral projective image $D$, then the maximal cyclic subgroup of $D$ is not a $p$-group,  
\item the map $G_{\Q} \rar \Aut{T_pE}$ given by the Galois action on the $p$-adic Tate module of $E$ is onto,  
\item $E$ is ordinary at $p$ and $\rho(\Fr_p)$ has distinct eigenvalues mod $\mathfrak{p}$.
\end{enumerate} 
If $L(E,\rho,1) \neq 0$, then the group $\mrm{Hom}_{\widehat{\OO_{\mfk{p}}}[\mrm{Gal(K/\Q)}]}(\rho,\mrm{Sel}_{p^{\infty}}(E/K)\otimes_{\Z_p} \widehat{\OO_{\mfk{p}}})$ is finite. 
}

We also show that the following weakening of Question \ref{loefflerqn} admits a positive answer.  

\propn[weakerqn]{If one replaces ``for all but finitely many $\mathfrak{p}$'' with ``for infinitely many $\mathfrak{p}$'' in Question \ref{loefflerqn}, then the answer is positive.}

However, we prove that the general answer to Question \ref{loefflerqn} is negative. Among other counter-examples, we find:

\propn[inf-cex]{(See Section \ref{sect-theoD}) Let $f \in \mathcal{S}_k(\Gamma_1(N))$ be a newform with character $\varepsilon_f$ without complex multiplication. Suppose that one of the following conditions holds:
\begin{itemize}[noitemsep,label=\tiny$\bullet$]
\item $k$ is even and $f$ has an inner twist by some odd primitive Dirichlet character. 
\item $k \geq 3$ is odd, and $f$ has an inner twist by two distinct non-trivial Dirichlet characters $\chi_1,\chi_2$ such that the restrictions to $\ker{\varepsilon_f}$ of $\chi_1$ and $\chi_2$ are both non-trivial and distinct. 
\end{itemize}
Then there exist a quadratic number field $K$, a primitive odd Dirichlet character $\varepsilon$, and an infinite collection $\mathcal{G}$ of weight one newforms $g$ with primitive character $\varepsilon$ satisfying the following properties:
\begin{itemize}[noitemsep,label=$-$] 
\item every $g \in \mathcal{G}$ has either real or complex multiplication by $K$,
\item for every distinct $g,g' \in \mathcal{G}$, $g$ and $g'$ are not quadratic twists of each other,
\item for every $g \in \mathcal{G}$, the answer to Question \ref{loefflerqn} for $(f,g)$ is negative.
\end{itemize}  
}

\propn[counterex]{(See Proposition \ref{special-counter} \ref{S3det-CM}) There exist newforms $f \in \mathcal{S}_2(\Gamma_0(63)), g \in \mathcal{S}_1(\Gamma_1(1452))$ with coefficients respectively in $\Q(\sqrt{3}),\Q(\sqrt{-3})$ satisfying the following property: for all but finitely many prime ideals $\mathfrak{p}$ of $L=\Q(e^{2i\pi/12})$ with residue characteristic $p \equiv 5,7 \pmod{12}$, there is no $\sigma \in G_{\Q(\mu_{p^{\infty}})}$ such that $\ker\left(\sigma-\mrm{id}\mid T_{f,\mathfrak{p}} \otimes_{\OO_{L_{\mathfrak{p}}}} T_{g,\mathfrak{p}}\right) \otimes_{\OO_{L_{\mathfrak{p}}}} L_{\mathfrak{p}}$ is a $L_{\mathfrak{p}}$-line.}

\rems{\begin{enumerate}[noitemsep,label=(\roman*)] 
\item The statement of Theorem \ref{twistedbsd} can likely be generalized by replacing elliptic curves with, for instance, abelian varieties attached to newforms of weight $2$, with only minimal changes to the proof. One would need to adapt the big image hypothesis, for instance by ensuring that the conclusion of Theorem \ref{weakerconds} holds. 
\item As explained in \cite[Remark 11.7.5(i)]{KLZ15}, applying Theorem \ref{twistedbsd} at a single prime shows that $\mrm{Hom}_{\OO[G_{\Q}]}(\rho,E(K) \otimes \OO)$ is zero. The point of Theorem \ref{twistedbsd} is to prove that, in addition, the $\rho$-part of $\Sha(E/K)[p^{\infty}]$ is finite for certain primes $p$. In fact, after removing finitely many primes by Serre's open image theorem \cite[Th\'eor\`eme 3]{Serre-image-ouverte}, Theorem \ref{twistedbsd} applies to any prime $p$ which is ordinary for $E$ and such that $\Fr_p$ does not lie in the kernel of the projective representation $\tilde{\rho}$ attached to $\rho$. The first condition excludes a set of primes of density zero by \cite[Cor. 2 au Th\'eor\`eme 20]{Serre-Cebotarev} (which is nonetheless infinite by \cite[Theorem 1]{Elkies}), and the second condition excludes, by Chebotarev's theorem, a set of primes of density $\frac{1}{|\mrm{im}(\tilde{\rho})|}$. In other words, Theorem B applies to a set of primes of density $1-\frac{1}{|\mrm{im}(\tilde{\rho})|}$, which is always at least $\frac{3}{4}$. 
\item Combining the first two points, when generalizing Theorem \ref{twistedbsd} to abelian varieties attached to weight two newforms, Theorem \ref{weakerqn} could imply that, provided that the special value of the $L$-function is not zero, the relevant Mordell-Weil group is finite.   
\item When $f \in \mathcal{S}_k(\Gamma_1(N))$, and $g \in \mathcal{S}_l(\Gamma_1(M))$ are newforms with characters $\varepsilon_f,\varepsilon_g$ respectively, the Dirichlet character $\varepsilon_f\varepsilon_g$ has the same parity as $k+l$. In particular, when $k$ is even and $l=1$, $\varepsilon_f\varepsilon_g$ is odd, so it is non-trivial. Thus, the counter-examples in Theorems \ref{inf-cex} and \ref{counterex} are not trivial.
\item We will see in Section \ref{sect-theoA} that if the answer to Question \ref{loefflerqn} is negative for a couple $(f,g)$, then it is also negative for every couple $(f,g')$, where $g'$ is any quadratic twist of $g$. This is why we require in Theorem \ref{inf-cex} that the weight one newforms should not be quadratic twists of each other.  
\item In the counter-example of Theorem \ref{counterex}, $g$ has complex multiplication, $f$ has even weight and an odd inner twist. While the answer to Question \ref{loefflerqn} is negative, this does not follow from the proof of Theorem \ref{inf-cex}. Further counter-examples are given in Section \ref{sect-counterex}. 
\end{enumerate}}

The failure of the hypothesis needed to apply the Euler system method in certain situations raises the following question. Recall that if $f$ is a newform with coefficients in some number field $L \subset \C$, and $\mfk{p} \subset \OO_L$ is a maximal ideal of residue characteristic $p$, we say (as in \cite[Definition 7.2.6]{KLZ15}) that $f$ is \emph{$\mfk{p}$-distinguished} if the semi-simplification of the restriction of $T_{f,\mfk{p}}/\mfk{p}$ to $G_{\Q_p}$ is the sum of two \emph{distinct} characters. For instance, in Theorem \ref{twistedbsd}, if $g$ denotes the weight one newform attached to $\rho$, then $p$ does not divide the level of $g$ by the first condition, and the second part of the fourth condition means precisely that $g$ is $\mfk{p}$-distinguished. 

Consider a couple $(f,g)$ which is a counter-example to Question \ref{loefflerqn}, for instance one of the couples given in Section \ref{sect-counterex}. Assume for the sake of simplicity that $f$ has weight two, and let $L_0 \subset \C$ be a number field containing all the Fourier coefficients of $f$ and all the Frobenius eigenvalues of $g$. We do not see any reason why the central archimedian Rankin--Selberg $L$-value $L(f \otimes g,1)$ should vanish, so let us assume that $L(f \otimes g,1) \neq 0$.

Let $\mfk{p} \subset L_0$ be a prime ideal with sufficiently large residue characteristic $p$: in this case, the $G_{\Q}$-modules $T_{f,\mfk{p}}/\mfk{p}$ and $T_{g,\mfk{p}}/\mfk{p}$ are absolutely irreducible (see Section \ref{sect-theoA}). Assume that $f$ is ordinary at $\mfk{p}$: then it is $\mfk{p}$-distinguished by \cite[Remark 7.2.7]{KLZ15}. Assume furthermore that $g$ is $\mfk{p}$-distinguished, and that there is no $\sigma \in G_{\Q(\mu_{p^{\infty}})}$ such that $\ker\left(\sigma-\mrm{id}\mid T_{f,\mathfrak{p}} \otimes_{\OO_{L_{\mathfrak{p}}}} T_{g,\mathfrak{p}}\right) \otimes_{\OO_{L_{\mathfrak{p}}}} L_{\mathfrak{p}}$ is a $L_{\mathfrak{p}}$-line\footnote{In the situation of Theorem \ref{counterex}, it can be shown that, for all primes $p \equiv 5,7 \pmod{12}$ such that either $p \equiv 5 \pmod{12}$ or $p$ is not a square modulo $11$, except a subset of density zero, there is a prime $\mfk{p}$ of $L_0=\Q(e^{\frac{2i\pi}{12}})$ satisfying the previous assumptions.}. 

Then should one expect, as in \cite[Theorem 11.7.3]{KLZ15}, that the Nekov\'a\u r cohomology group $\tilde{H}^2\left(\OO_L\left[\frac{1}{N_fN_gp}\right],T_{f,\mfk{p}} \otimes T_{g,\mfk{p}};\Delta^{BK}\right)$ is finite\footnote{The statement of the theorem in \emph{loc. cit.} seems to contain a slight typo; one should (probably) read $M_{\OO}(f \otimes g)^{\ast}$ instead of $M_{\OO}(f \otimes g)$.}? How to go about proving such a claim without the theory of Euler systems does not seem clear at all.  \\

%


\textbf{Acknowledgements} \, This question arose in the context of my PhD thesis, and I am grateful to my supervisor Lo\"ic Merel for encouraging me to work on this project. His careful reading and many comments helped me vastly improve the original version of this text. I would like to thank Pierre Parent and Jan Vonk for their comments as well. I would also like to thank David Loeffler for his comments and his help regarding the application of the results of this paper to the twisted Birch and Swinnerton-Dyer conjecture. I would finally like to thank the anonymous reviewers for their attentive reading of this article and their many remarks, comments and suggestions. This research was partially funded by the European Union, ERC Starting Grant 101076941 'GAGARIN'.

\section{The big image theorem for one modular form}
\label{sect-setup}

We first recall well-known results on twists and Galois conjugates of modular forms. This lets us define inner twists, which are crucial to stating the large image theorem for one modular form.  

The results on twists are mostly taken from \cite{AL78}, the lemmas on inner twists come from \cite{Momose}, and the presentation of the results on Galois images of this section comes almost entirely from \cite[Section 2]{bigimage}.  

We start with twists of newforms by Dirichlet characters. Although the author is not aware of a complete proof of the result below exactly as we state it, it is an easy consequence of the classical results of \cite[\S 3]{AL78} (on twists of modular forms by Dirichlet characters) and \cite[\S 4.6]{Miyake} (about the properties of newforms).

\textbf{Notation} \; If $q$ is a prime number and $N$ is a nonzero integer, $v_q(N)$ is the $q$-adic valuation of $N$, i.e. the greatest integer $r \geq 0$ such that $q^r \mid N$. \\

\prop[twist-exists]{Let $f \in \mathcal{S}_k(\Gamma_1(N))$ be a newform with character $\varepsilon_f$ of conductor $M$. Let $\chi$ be a primitive Dirichlet character of conductor $C$. There exists a unique new, cuspidal eigenform $f \otimes \chi \in \mathcal{S}_k(\Gamma_1(L))$ such that, for every prime number $p\nmid NC$, $a_p(f \otimes \chi)=a_p(f)\chi(p)$. We say that $f\otimes \chi$ is the \emph{twist} of $f$ by $\chi$. Moreover: 
\begin{enumerate}[label=$($\alph*$)$,noitemsep]
\item \label{twist-character} The character of $f \otimes \chi$ is the primitive Dirichlet character associated to $\varepsilon_f\chi^2$.
\item \label{level-twist-basic} $L$ divides the least common multiple of $N,CM,C^2$.
\item \label{level-twist-coprime} For every prime number $q \nmid C$, one has $v_q(L)=v_q(N)$.
\item \label{level-twist-local} Let $q$ be a prime number, and write $\chi=\chi_q\chi'$, where the conductor of $\chi_q$ (resp. $\chi'$) is a power of $q$ (resp. coprime to $q$). Then $v_q(L)$ only depends on $f$ and $\chi_q$.
\end{enumerate}}

\rem{If $\psi$ comes from the primitive Dirichlet character $\psi_0$, we will also denote by $f \otimes \psi$ the newform $f \otimes \psi_0$. }

\prop{Let $f \in \mathcal{S}_k(\Gamma_1(N))$ and write its $q$-expansion $f(\tau)=\sum_{n \geq 1}{a_n(f)q^n}$. Let $\sigma$ be any automorphism of $\C$. Then $\sigma(f)(\tau) := \sum_{n \geq 1}{\sigma(a_n(f))q^n} \in \mathcal{S}_k(\Gamma_1(N))$. Moreover:
\begin{enumerate}[label=$($\alph*$)$,noitemsep]
\item If $f$ has character $\varepsilon_f$, then $\sigma(f)$ has character $\sigma\circ \varepsilon_f$. 
\item If $f$ is a newform, $\sigma(f)$ is a newform. 
\item If $f$ is a newform, there is a number field containing all its coefficients. 
\end{enumerate} }

\bigskip

Let $f \in \mathcal{S}_k(\Gamma_1(N))$ be a newform with (primitive) character $\varepsilon_f$. Let $f(\tau) = \sum_{n \geq 1}{a_n(f)q^n}$ be its $q$-expansion and let $L=\Q(a_n(f)\mid n \geq 1) \subset \C$. \\

\defi{Let $\gamma: L \rar \C$ be a field homomorphism and $\chi$ be a primitive Dirichlet character. We say that $f$ has an \emph{inner twist} by $(\gamma,\chi)$ if the newforms $\gamma(f)$ and $f \otimes \chi$ are equal. }

\lem[momose15]{If $f$ has an inner twist by $(\gamma,\chi)$, then one can write $\chi=\lambda\varepsilon_f^t$, where $t \in \Z$ and $\lambda$ is quadratic. Moreover, the conductor of $\chi$ divides $N$ and one has $\gamma(L) = L$.}

\demo{Apart from the bound on the conductor of $\chi$, this is \cite[Lemma 1.5 (i)]{Momose}. 

Let $C$ be the conductor of $\chi$, and let $p$ be a prime such that $v_p(C) > v_p(N)$. Then by \cite[Theorem 4.1]{AL78} and Proposition \ref{twist-exists}\ref{level-twist-local}, the $p$-adic valuation of the conductor of $f \otimes \chi$ is $v_p(N)$, and is also $2v_p(C)$, and we get a contradiction.   
}

\medskip
\noindent
\textbf{Notation} \; Let $\Gamma_f$ denote the set of $(\gamma,\chi)$ by which $f$ has an inner twist. 

\medskip

\prop[abelian-dicho]{The formula $(\alpha,\chi) \cdot (\beta,\psi) = (\alpha\circ \beta, \chi\alpha(\psi))$ endows $\Gamma_f$ with an abelian group law. Moreover, the group homomorphism \[\pi: (\gamma,\chi) \in \Gamma_f \longmapsto \gamma \in \Aut{L}\] satisfies exactly one of the following properties:
\begin{itemize}[label=\tiny$\bullet$]
\item $\pi$ is injective with abelian image. Let $F$ be the subfield of $L$ fixed by $\mrm{im}(\pi)$ (which we identify to $\Gamma_f$): the extension $L/F$ is abelian with Galois group $\mrm{im}(\pi) \simeq \Gamma_f$. 
\item There exists some $(\mrm{id},\eta) \in \Gamma_f$ where $\eta$ is the character of a quadratic field $K$ which is only ramified at primes dividing $N$. If $k \geq 2$, then $K$ is imaginary. 
\end{itemize}}

\demo{Checking that we have indeed defined a group law on $\Gamma_f$ such that $\pi$ is a homomorphism is formal. 

Let $(\sigma,\alpha),(\tau,\beta) \in \Gamma_f$. Then we can write $\alpha=\alpha_2\varepsilon_f^a$, $\beta=\beta_2\varepsilon_f^b$, where $a,b \in \Z$ and $\alpha_2,\beta_2$ are quadratic Dirichlet characters. Then $(\tau \circ \alpha)\cdot \alpha^{-1}=(\tau \circ \varepsilon_f)^a \cdot (\varepsilon_f^{-1})^a=(\beta^2)^a=\varepsilon_f^{2ab}$: by symmetry, $(\tau\circ \alpha)\cdot \beta=\alpha\cdot (\sigma\circ\beta)$. %

Let $\gamma$ be the primitive Dirichlet character associated to $(\tau\circ\alpha)\cdot \beta$, so that 
\[(\sigma,\alpha) \cdot (\tau,\beta)=(\sigma\circ\tau,\gamma),\quad (\tau,\beta)\cdot (\sigma,\alpha) = (\tau\circ \sigma,\gamma).\] 

The newforms $\sigma(\tau(f)),\tau(\sigma(f))$ are equal to the newform $f \otimes \gamma$, so they are equal. Hence the automorphisms $\sigma\circ\tau$ and $\tau\circ\sigma$ of $L$ agree on every Fourier coefficient of $f$, so $\sigma\circ\tau=\tau\circ\sigma$. Therefore, $(\sigma,\alpha) \cdot (\tau,\beta)= (\tau,\beta)\cdot (\sigma,\alpha)$, i.e. $\Gamma_f$ is abelian. 

If $\pi$ is injective, its image is an abelian subgroup of $\Aut{L}$. Let $F$ be the field fixed by $\mrm{im}(\pi)$, then it is well-known (e.g. \cite[VI, Theorem 1.8]{LangAlg}) that $L/F$ is Galois with Galois group $\mrm{im}(\pi) \simeq \Gamma_f$ which is abelian. 

Suppose that $\pi$ is not injective. Then there exists a non-trivial Dirichlet character $\chi$ such that $f \otimes \chi=f$. Thus $\varepsilon_f\chi^2=\varepsilon_f$, so $\chi$ is quadratic. The rest follows from \cite[(4.4, 4.5)]{Antwerp5-Ribet} in weight $k \geq 2$. 
}

In the latter case, we say that $f$ has \emph{complex} (resp. \emph{real}) \emph{multiplication by $K$} if $K$ is imaginary (resp. real). This is often shortened as ``CM'' (resp. ``RM''), and we will omit $K$ from the notation if it does not matter. \\

\noindent
\textbf{Notation} \; If $f$ does not have CM (or RM), let $F$ be the subfield of $L$ fixed by the image of $\pi$, and let $H/\Q$ be the abelian field extension such that $G_H = \bigcap_{(\gamma,\chi) \in \Gamma_f}{\ker{\chi}}$, where $\chi$ is viewed as a character of $G_{\Q}$ through the cyclotomic character $G_{\Q} \rar \hat{\Z}^{\times}$. To be completely explicit, we identify \emph{any} primitive Dirichlet character $\eta$ of level $D$ with the character $G_{\Q} \rar \C^{\times}$ unramified outside $D$ mapping, for all primes $p\nmid D$, the arithmetic Frobenius at $p$ to $\eta(p)$. 

\rems{\begin{itemize}[noitemsep,label=\tiny$\bullet$]
\item If $c$ denotes the complex conjugation, one always has $(c,\varepsilon_f^{-1}) \in \Gamma_f$ (see \cite[p. 21]{Antwerp5-Ribet}). As a consequence, the Galois character $\varepsilon_f: G_{\Q} \rar L^{\times}$ factors through $\mrm{Gal}(H/\Q)$. Furthermore, the number field $F$ is fixed by the complex conjugation, so, since $L$ is totally real or CM by \cite[Proposition 3.2]{Antwerp5-Ribet}, $F$ is totally real.  
\item When $f$ is replaced with a quadratic twist, it is elementary to check that $L$ and $\Gamma_f$ are unchanged, so the fields $H$ and $F$ remain the same as well. 
\end{itemize}
}

\medskip
\noindent

Let $\mathfrak{p} \subset \OO_L$ be a prime ideal with finite residue characteristic $p$. It is well-known (see \cite[Theorems 2.1, 2.3]{Antwerp5-Ribet}) that there exists a group homomorphism $\rho_{f,\mathfrak{p}}: G_{\Q} \rar \GL{\OO_{L_{\mathfrak{p}}}}$ satisfying the following properties:

\begin{itemize}[noitemsep,label=\tiny$\bullet$]
\item $\rho_{f,\mathfrak{p}}$ is continuous, irreducible, and unramified away from $pN$.
\item $\det{\rho_{f,\mathfrak{p}}}=\varepsilon_f \omega_p^{k-1}$, where $\omega_p$ is the $p$-adic cyclotomic character. 
\item For every prime $\ell \nmid pN$, $\Tr{\rho_{f,\mathfrak{p}}(\Fr_{\ell})} = a_{\ell}(f)$, where $\Fr_{\ell}$ is the arithmetic Frobenius.  
\end{itemize}

\rem{The $\rho_{f,\mfk{p}}$ that we consider is the one considered throughout the paper \cite{Antwerp5-Ribet}. It is, however, \emph{dual} to the one considered in \cite[Definition 2.1.1]{bigimage}. As pointed out in \cite[Remark 2.1.2]{bigimage}, this is not important for our purposes since we are only considering the image of the representation. }

\prop{(Momose \cite[Theorem 4.2]{Momose}, Ribet \cite[\S 3]{Ribet2}, Ghate--Gonz\'alez-Jim\'enez--Quer \cite[Corollary 4.7]{GGJQ}) 

Assume that $k \geq 2$, that $p \nmid 2N$ is unramified in $L$ and $f$ does not have CM. Then $\rho_{f,\mathfrak{p}}(G_H)$ has a $\GL{\OO_{L_{\mathfrak{p}}}}$-conjugate contained in $\{M \in \GL{\OO_{F_{\mathfrak{p}}}},\, \det{M} \in \Z_p^{\times}\}$. 

Furthermore, for all but finitely many $\mathfrak{p}$, after conjugating $\rho_{f,\mathfrak{p}}$ by a suitable element of $\GL{\OO_{L_{\mathfrak{p}}}}$, the group homomorphism $(\rho_{f,\mathfrak{p}},\omega_p): G_H \rar \GL{\OO_{F_{\mathfrak{p}}}} \times_{\OO_{F_{\mathfrak{p}}}^{\times}} \Z_p^{\times}$ is onto (where the left map in the fibre product is the determinant, and the right one is $x \in \Z_p^{\times} \longmapsto x^{k-1} \in \OO_{F_{\mathfrak{p}}}^{\times}$). }

\demo{This is \cite[Theorem 2.2.2]{bigimage}. }

\rem{While we will not be considering it in the rest of the paper, let us briefly discuss the situation when $k \geq 2$ and $f$ has CM by an imaginary quadratic field $K$. The following claims can essentially be deduced from the contents of \cite[\S 4]{Antwerp5-Ribet} (using Serre's theorem on rational abelian semi-simple representations \cite[Chap. III, \S 3.1, Theorem and Remark 1]{McGill}), and in particular the fact that $f$ is constructed from a Gr\"ossencharakter $\psi$ as in \cite[pp. 34--35]{Antwerp5-Ribet}. 

Let $\mfk{p}$ be a maximal ideal of $\OO_L$ with residue characteristic $p$. Let $m \geq 1$ be the norm of the modulus of definition of $\psi$, and $H$ be the ray class field of $K$ modulo $m$. Assume for simplicity that $p$ is prime to $6m$ and unramified in $K/\Q$. Then $H/\Q$ is unramified above $p$, $(\rho_{f,\mfk{p}})_{|G_K}$ is abelian, and one has $\rho_{f,\mfk{p}}(I_{K,p})=\rho_{f,\mfk{p}}(G_H)$, where $I_{K,p}$ is the subgroup of $G_K^{ab}$ generated by the inertia subgroups of prime ideals of $K$ above $p$. 

Moreover, let $\mfk{a}: \mathbb{A}_K^{\times}/K^{\times} \rar G_K^{ab}$ denote the class field theory reciprocity homomorphism (mapping a uniformizer to an arithmetic Frobenius), then one has $\mfk{a}((\OO_K \otimes \Z_p)^{\times})=I_{K,p}$ and $\left(\rho_{f,\mfk{p}} \circ \mfk{a}\right)_{|(\OO_K \otimes \Z_p)^{\times}}$ is conjugate to the $(1-k)$-th power of the representation \[(\OO_K \otimes \Z_p)^{\times} \subset \Auti{\Z_p}{\OO_K \otimes \Z_p} \simeq \GL{\Z_p} \subset \GL{\OO_{L_{\mfk{p}}}}.\] }

\section{Special elements in the image of the convolution of two modular forms}
\label{sect-theoA}

Let $p$ be an odd prime. Let $\OO$ be the ring of integers of a finite extension $K$ of $\Q_p$, with uniformizer $\varpi$ and residue field $k$. Let $T$ be a free $\OO$-module of finite rank, endowed with a continuous and $\OO$-linear action of $G_{\Q}$. Let $G$ be the image in $\Aut{T}$ of $G_{\Q(\mu_{p^{\infty}})}$. Consider the following statements:  

\begin{enumerate}[noitemsep]
\item[(N)] $-\mrm{id} \in G$. 
\item[(gI)] $T \otimes_{\OO} K$ is an irreducible $K[G]$-module. 
\item[(rI)] $T \otimes_{\OO} k$ is an irreducible $k[G]$-module.
\item[(wE)] There exists $u \in G$ such that $\ker(u-\mrm{id}) \otimes K$ has dimension one. 
\item[(sE)] There exists $u \in G$ such that $T/(u-1)T$ is free of rank one. 
\end{enumerate}

\defi{We say that $T$ has Galois image \emph{of Euler type} if it satisfies (N), (gI) and (wE).\\ $T$ has  \emph{Euler-adapted} Galois image if it satisfies (N), (rI) and (sE).}

Note that (rI) implies (gI) (since any reducible $K[G]$-submodule can be reduced in $k$), and that (sE) implies (wE). 

The statements (gI) and (wE) (resp. (rI) and (sE)) are the two parts of $\mrm{Hyp}(\Q(\mu_{p^{\infty}}), T \otimes K)$ (resp. $\mrm{Hyp}(\Q(\mu_{p^{\infty}}), T)$) in the notation of \cite[Section 4.1]{bigimage}. The condition (N) is not mentioned explicitly in \cite{bigimage}, but is used in applications of this article to Euler systems, such as \cite[Section 11]{KLZ15}. Since it is not difficult to study, we include it here.   \\

Let $k \geq 2$. Let $f \in \mathcal{S}_k(\Gamma_1(N_f))$, $g \in \mathcal{S}_1(\Gamma_1(N_g))$ be newforms with characters $\varepsilon_f$ and $\varepsilon_g$ such that $f$ has no complex multiplication. Let $L,F,H$ be the number fields associated to $f$ in the previous section and $L_0 \subset \C$ be a number field containing $L$ and all the Fourier coefficients of $g$.  

If $\mathfrak{p} \subset \OO_{L_0}$ is a maximal ideal, we can construct the representation $T_{f,\mathfrak{p}}$ of $G_{\Q}$ as in the previous section by extension of scalars: it is a free $\OO_{(L_0)_{\mathfrak{p}}}$-module of rank $2$ with a continuous linear action of $G_{\Q}$. When the residue characteristic of $\mathfrak{p}$ does not divide $2N_f$ or the discriminant of $L$, we endow $T_{f,\mathfrak{p}}$ with a basis such that the image of $G_H$ is contained in $\GL{\OO_{F_{\mathfrak{p}}}}$. 

By the work of Deligne and Serre \cite[Th\'eor\`emes 4.1, 4.6]{Del-Ser} and \cite[Lemma 5.3]{Antwerp5-Ribet}, there is an odd irreducible Artin representation $\rho_g: G_{\Q} \rar \GL{\OO_{L_0}}$ attached to $g$. Let $\tilde{\rho}_g: G_{\Q} \rar \PGL{L_0}$ denote the corresponding projective representation and let $Q=|\tilde{\rho}_g(G_{\Q})|$. Let $T_g$ be a free $\OO_{L_0}$-module of rank $2$ where $G_{\Q}$ acts as $\rho_g$. Let $T_{f,g,\mathfrak{p}} = T_{f,\mathfrak{p}} \otimes_{\OO_{L_0}} T_{g}$. After enlarging $L_0$, we assume that $L_0$ contains the eigenvalues of every matrix in $\rho_g(G_{\Q})$. 

\rem{The subgroup $\tilde{\rho}_g(G_{\Q})$ is a finite subgroup of $\PGL{L_0}$. Such subgroups admit a well-known classification (see e.g. \cite[Proposition 16]{Serre-image-ouverte}) which we now recall, since it is used throughout the paper. More generally, if $K$ is a characteristic zero field, any finite subgroup of $\PGL{K}$ satisfies exactly one of the following conditions:
\begin{itemize}[noitemsep,label=\tiny$\bullet$]
\item it is cyclic, 
\item it is dihedral of order $2n$ for some $n \geq 2$ (we call ``dihedral'' the abelian group $\F_2^{\oplus 2}$), 
\item it is isomorphic to one of the classical groups $A_4$, $S_4$, or $A_5$.
\end{itemize}
The group $\tilde{\rho}_g(G_{\Q})$ is not cyclic, because otherwise the quotient of $\rho_g(G_{\Q})$ by its center is cyclic, so $\rho_g(G_{\Q})$ is abelian, which contradicts the irreducibility of $\rho_g$. 
}

\medskip
\noindent
\defi{A maximal ideal $\mathfrak{p} \subset \OO_{L_0}$ is \emph{good} for $(f,g)$ if the following conditions are satisfied:
\begin{itemize}[noitemsep,label=\tiny$\bullet$]
\item Its residue characteristic $p$ does not divide $30N_fN_g$ or ramify in $L$. 
\item $(T_{f,\mathfrak{p}},\omega_p): G_H \rar \GL{\OO_{F_{\mathfrak{p}}}} \times_{\OO_{F_{\mathfrak{p}}}^{\times}} \Z_p^{\times}$ is surjective, where the fibre product maps are the determinant and the $(k-1)$-th power.
\end{itemize}
Note that, by the results of the previous section, all but finitely many maximal ideals of $\OO_{L_0}$ are good. }

\prop[jointimage-H]{Let $\mathfrak{p}$ be a good prime with residue characteristic $p$. The image of $G_{H(\mu_{p^{\infty}})}$ in $\Aut{T_{f,\mathfrak{p}}} \times \Aut{T_g}$ is $\SL{\OO_{F_{\mathfrak{p}}}} \times \rho_g(G_H)$. Moreover, for any $\sigma \in G_{\Q}$, there exists $\sigma' \in G_{\Q(\mu_{p^{\infty}})}$ such that $\rho_g(\sigma)=\rho_g(\sigma')$ and $\sigma,\sigma'$ have the same image in $\mrm{Gal}(H/\Q)$. }

\demo{Let $\OO$ denote the ring of integers of $F_{\mfk{p}}$. Let $I$ denote the image of $G_{H(\mu_{p^{\infty}})}$ in $\Aut{T_{f,\mathfrak{p}}} \times \Aut{T_g}$. It is clear that $I \subset \SL{\OO} \times \rho_g(G_H)$. 

Let $p_1: I \rar \SL{\OO},\, p_2: I \rar \rho_g(G_H)$ be the two projections. Because $\mathfrak{p}$ is a good prime, $p_1$ is surjective. 

We claim that for any finite Galois extension $M/\Q$ unramified at $p$, $G_{\Q(\mu_{p^{\infty}})} \rar \mrm{Gal}(M/\Q)$ is surjective. 
Indeed, it is enough to show that $\Q(\mu_{p^{\infty}}) \cap M = \Q$. Now, $\Q(\mu_{p^{\infty}}) \cap M$ is a finite extension of $\Q$, so it is equal to $\Q(\mu_{p^n}) \cap M$ for some $n \geq 1$. By assumption, $M/\Q$ is unramified at $p$, so $\Q(\mu_{p^n}) \cap M$ is unramified at $p$ over $\Q$. Since $\Q(\mu_{p^n})/\Q$ is totally ramified at $p$, so is the intermediate extension $(\Q(\mu_{p^n}) \cap M)/\Q$. Thus, $\Q(\mu_{p^{n}}) \cap M$ is both unramified and totally ramified over $\Q$ at $p$, so $\Q=\Q(\mu_{p^{n}}) \cap M=\Q(\mu_{p^{\infty}}) \cap M$. 

The extension $M/\Q$ generated by $H$ and the splitting field of $\rho_g$ is unramified at $p$ since $p \nmid N_fN_g$, so $G_{\Q(\mu_{p^{\infty}})} \rar \mrm{Gal}(M/\Q)$ is surjective. Hence, for any $\sigma \in G_{\Q}$, there is $\sigma' \in G_{\Q(\mu_{p^{\infty}})}$ such that $\rho_g(\sigma)=\rho_g(\sigma')$ and $\sigma' \in \sigma G_H$. In particular, $p_2$ is onto and $\rho_g(G_{H'(\mu_{p^{\infty}})})=\rho_g(G_{H'})$ for any subextension $\Q \subset H' \subset H$. 

Therefore, it only remains to show that $I$ contains $\{1\} \times \rho_g(G_H)$. The group $p_2(I)=\rho_g(G_H)$ contains the normal subgroup $G_1=p_2(\ker{p_1})$. Let $\pi: \SL{\OO} \rar \rho_g(G_H)/G_1$ be the continuous surjective homomorphism defined as follows: for each $x \in \SL{\OO}$, pick $y \in \rho_g(G_H)$ such that $(x,y) \in I$, and let $\pi(x)=y \pmod{G_1}$. This is well-defined because any two choices of $y$ differ by an element of $G_1$. We need to show that $G_1=\rho_g(G_H)$.

Recall that, for a Hausdorff topological group $G$, its \emph{derived subgroup} is the closed subgroup $G'$ generated by its commutators, i.e. the elements of the form $ghg^{-1}h^{-1}$ with $g,h \in G$. The subgroup $G'$ is preserved by every continuous automorphism of $G$, and it is contained in the kernel of every continuous map $G \rar A$ where $A$ is a Hausdorff topological abelian group. 

For $u \in \OO$, let $M_u := \begin{pmatrix}1 & u\\0 & 1\end{pmatrix} \in \SL{\OO}$ and $M'_u=M_u^T \in \SL{\OO}$. By \cite[(1.2.11)]{HOM}, the matrices $M_{u}, M'_u$ for $u \in \OO$ generate $\SL{\OO}$. Furthermore, if $W=\begin{pmatrix}0 & -1\\1 & 0\end{pmatrix}$, one checks that $M'_u = WM_uW^{-1}$, so the normal subgroup of $\SL{\OO}$ generated by the $M_u$ is $\SL{\OO}$. 

Since $p > 3$, we can find $\lambda \in \OO$ such that $\lambda(\lambda^2-1)$ is a unit. For any $u \in \OO$, one has \[M_u := \begin{pmatrix}1 & u\\0 & 1\end{pmatrix} = \begin{pmatrix}\lambda & 0\\0 & \lambda^{-1} \end{pmatrix}\begin{pmatrix}1 & \frac{u}{\lambda^2-1}\\0 & 1\end{pmatrix}\begin{pmatrix}\lambda^{-1} & 0\\0 & \lambda\end{pmatrix}\begin{pmatrix}1 & \frac{-u}{\lambda^2-1}\\0 & 1\end{pmatrix} \in \SL{\OO}'.\] 

Thus $\SL{\OO}'$ contains the normal subgroup generated by the $M_u$, so $\SL{\OO}'=\SL{\OO}$. In particular, every continuous group homomorphism from $\SL{\OO}$ to a finite abelian group is trivial. 

Moreover, if $G$ is a finite group with cardinality $d$ prime to $p$ and $\varphi: \SL{\OO} \rar G$ is a group homomorphism, then, for any $u \in \OO$, one has $\varphi(M_{u})=\varphi(M_{u/d}^d)=\varphi(M_{u/d})^d=1_G$. So $\ker{\varphi}$ contains the normal subgroup generated by all the $M_u$, so $\ker{\varphi}=\SL{\OO}$ i.e. $\varphi$ is trivial. 

Therefore, for any finite group $G$ whose composition factors are either cyclic or simple groups of cardinality prime to $p$, any continuous homomorphism $\varphi: \SL{\OO} \rar G$ is trivial. Since $\tilde{\rho}_g(G_H)$ is cyclic, dihedral, or isomorphic to $A_4, S_4$ or $A_5$, the composition factors of $\rho_g(G_H)$ are cyclic or isomorphic to $A_5$, so the same holds for $\rho_g(G_H)/G_1$. 
Since $p$ does not divide $|A_5|=60$, $\pi$ is therefore trivial and surjective, so $G_1=\rho_g(G_H)$. }

\cor[papier]{Let $\mathfrak{p}$ be a good prime with residue characteristic $p$ and $\sigma \in G_{\Q}$. There exists $\alpha \in \OO_{L,\mathfrak{p}}^{\times}$ such that for all $(\gamma,\chi) \in \Gamma_f$, one has $\gamma(\alpha)=\chi(\sigma)\alpha$. If moreover $\sigma \in G_{\Q(\mu_{p^{\infty}})}$, the image of $\sigma G_{H(\mu_{p^{\infty}})}$ in $\Aut{T_{f,\mathfrak{p}}} \times \Aut{T_g}$ is $\left[\begin{pmatrix}\alpha & 0\\0 & \alpha^{-1}\varepsilon_f(\sigma)\end{pmatrix}\SL{\OO_{F_{\mathfrak{p}}}}\right] \times \left[\rho_g(\sigma)\rho_g(G_H)\right]$.}

\demo{Let us briefly recall why $\alpha$ exists. For any $\gamma \in \mrm{Gal}(L/F)$, there is a unique primitive Dirichlet character $\chi_{\gamma}$ such that $(\gamma,\chi_{\gamma}) \in \Gamma_f$. The map $\gamma \in \mrm{Gal}(L/F) \longmapsto \chi_{\gamma}(\sigma) \in L^{\times}$ is a cocycle: by Hilbert 90 (for instance \cite[\S X.1, Proposition 2]{Serre-local}), it can be written as $\gamma \longmapsto \gamma(\beta)\beta^{-1}$ for some $\beta \in L^{\times}$. Since $p$ is unramified in $L$, $p$ generates the maximal ideal $\mfk{p}\OO_{L,\mfk{p}}$, so there is some $t \in \Z$ such that $\alpha=p^t\beta$ lies in $\OO_{L,\mathfrak{p}}^{\times}$ and satisfies the required conditions. 

By Proposition \ref{jointimage-H}, the image in $\Aut{T_{f,\mathfrak{p}}} \times \Aut{T_g}$ of $\sigma G_{H(\mu_{p^{\infty}})}$ is $I_{\sigma} \times [\sigma \rho_g(G_H)]$, where $I_{\sigma}$ is the image in $\Aut{T_{f,\mathfrak{p}}}$ of $\sigma G_{H(\mu_{p^{\infty}})}$. Now, $I_{\sigma} = \begin{pmatrix}\alpha & 0\\0 & \alpha^{-1}\varepsilon_f(\sigma)\end{pmatrix}\SL{\OO_{F_{\mathfrak{p}}}}$ by \cite[Corollary 2.2.3]{bigimage}.}

\rem{For every $(\gamma,\chi) \in \Gamma_f$, one has $\gamma\left(\frac{\alpha^2}{\varepsilon_f(\sigma)}\right)=\frac{\chi(\sigma)^2\alpha^2}{(\chi^2\varepsilon_f)(\sigma)}=\frac{\alpha^2}{\varepsilon_f(\sigma)}$, so $\frac{\alpha^2}{\varepsilon_f(\sigma)} \in \OO_{F,\mfk{p}}^{\times}$. In particular, one has \[\begin{pmatrix}\alpha & 0\\0 & \alpha^{-1}\varepsilon_f(\sigma)\end{pmatrix}\SL{\OO_{F_{\mathfrak{p}}}}=\left\{\alpha M\mid\; M \in \GL{\OO_{F_{\mfk{p}}}}, \det{M}=\frac{\varepsilon_f(\sigma)}{\alpha^2}\right\}.\]

It is clear from the proof that $\alpha$ is defined up to a $\mfk{p}$-adic unit of $F$ (i.e. an element of $\OO_{F,\mfk{p}}^{\times}$), and only depends on the image of $\sigma$ in $\mrm{Gal}(H/\Q)$.
}

\lem[irred-product-abstract]{Let $V,W$ be finite-dimensional vector spaces over a field $k$, and $G_1,G_2$ be subgroups of $\mrm{GL}(V)$ and $\mrm{GL}(W)$. Let $\Gamma$ be a subgroup of $\mrm{GL}(V) \times \mrm{GL}(W)$, containing $G_1 \times \{\mrm{id}\}$ and such that its second projection is $G_2$. Suppose that $V$ is an absolutely irreducible $k[G_1]$-module and $W$ is an irreducible $k[G_2]$-module. Then $V \otimes W$ is an irreducible $k[\Gamma]$-module. }

\demo{Let $U \subset V \otimes W$ be a proper $k[\Gamma]$-submodule. Let $W^{\ast}$ denote the dual vector space to $W$; if we let $G_1$ act trivially on $W$, by Schur's lemma, the evaluation map $\iota: W^{\ast} \rar \mathrm{Hom}_{G_1}(V \otimes W, V)$ is an isomorphism. Let $W'$ be the kernel of $w' \in W^{\ast} \mapsto \iota(w')_{|U}$, i.e. the set of linear forms $w': W \rar k$ such that $U \subset V \otimes (\ker{w'})$. By definition, $W'$ is a vector subspace of $W^{\ast}$. 

If $g \in G_2$ and $w' \in W'$, then there exists $\gamma=(g_1,g_2) \in \Gamma$ with $g_1 \in \GL{V}$. Thus one has $\iota(w' \circ g_2)=g_1^{-1}\circ \iota(w') \circ \gamma$. Since $U$ is stable under $\Gamma$, $\iota(w' \circ g_2)_{|U}$ is zero if and only if $g_1^{-1}\circ \iota(w')_{|U}$ is zero, if and only if $\iota(w')$ is zero: hence $W'$ is stable under $G_2$. Since $W$ is an irreducible $k[G_2]$-module, $W^{\ast}$ is irreducible under the action of $G_2$ as well, so $W'$ is either zero or $W^{\ast}$. 

The irreducible Jordan--H\"older components (cf. e.g. \cite[(13.7)]{CuRe}) of the $k[G_1]$-module $V \otimes W$ are isomorphic to $V$, so the same holds for $U$. Now, $V \otimes W \simeq V^{\oplus \dim{W}}$ is a completely reducible $k[G_1]$-module, so $U$ is a direct sum of irreducible $k[G_1]$-submodules by \cite[(15.2), Theorem 15.3]{CuRe}. Therefore $U \simeq V^{\oplus d}$ as a $k[G_1]$-module, with $d=\frac{\dim_k{U}}{\dim_k{V}}$. Since $U$ is properly contained in $V \otimes W$, one has $d < \dim_k{W}$, so that by Schur's lemma $d=\dim_k\,\mrm{Hom}_{k[G_1]}(U,V) < \dim_k{W^{\ast}}$. Hence $w' \in W^{\ast} \mapsto \iota(w')_{|U}$ is not injective, thus $W'$ is nonzero. Therefore $W'=W^{\ast}$, thus $U=0$ and the conclusion follows. 

}

\prop[only-we-se]{If $\mathfrak{p}$ is a good prime of residue characteristic $p$, $T_{f,g,\mathfrak{p}}$ satisfies conditions (N) and (gI). Moreover, if $Q/2$ is not a power of $p$, then condition (rI) holds.}

\demo{By Proposition \ref{jointimage-H}, $G_{H(\mu_{p^{\infty}})}$ contains an element $\sigma$ acting by $-I_2$ on $T_{f,\mathfrak{p}}$ and $I_2$ on $T_g$, which proves condition (N). 

The image $I$ of the action of $G_{\Q(\mu_{p^{\infty}})}$ in $\Aut{T_{f,\mfk{p}}} \times \Aut{T_g}$ contains $\SL{\OO_{F_{\mfk{p}}}} \times \{\mrm{id}\}$ and its second projection is $\rho_g(G_{\Q})$ by Proposition \ref{jointimage-H}.   

Since $\SL{\OO_{F_{\mfk{p}}}}$ acts absolutely irreducibly on $F_{\mfk{p}}^{\oplus 2}$, (gI) holds by Lemma \ref{irred-product-abstract}. Moreover, $\SL{\Z_p}$ acts absolutely irreducibly on $\F_p^{\oplus 2}$. Thus, if $\rho_g(G_{\Q})$ acts absolutely irreducibly on $T_g/\mathfrak{p}$, then Lemma \ref{irred-product-abstract} also implies condition (rI). 

Let $k$ be a finite extension of $\OO_{(L_0)_{\mfk{p}}}/\mfk{p}$ and assume that $T_g/\mfk{p} \otimes k$ contains a non-trivial invariant $k[G_{\Q}]$-submodule $U$: then $U$ is a line over $k$. Let $c \in G_{\Q}$ be the complex conjugation: the line $U$ is stable under $c$, so $U$ is equal to one of the eigenspaces of $c$. Since the characteristic polynomial of $\rho_g(c)$ is $(x+1)(x-1)$ and $2 \notin \mfk{p}$, there is an eigenspace $U_0$ for the action of $c$ on $T_g/\mfk{p}$ such that $U=U_0 \otimes k$. Then $U_0$ is stable under $G_{\Q}$, and $T_g/\mfk{p}$ is reducible under the action of $G_{\Q}$.  

Therefore, (rI) holds if $G_{\Q}$ acts irreducibly on $T_g/\mathfrak{p}$. 

Suppose that $G_{\Q}$ acts reducibly on $T_g/\mathfrak{p}$. This action factors through $\rho_g(G_{\Q})$, and its semisimplification is the direct sum of characters $\alpha, \beta: \rho_g(G_{\Q}) \rar (\OO_{L_0}/\mfk{p})^{\times}$. Thus $\gamma = \alpha/\beta$ is a character of $\rho_g(G_{\Q})$ and its kernel $K$ consists of matrices $g$ whose reduction modulo $\mathfrak{p}$ is conjugate to $\begin{pmatrix} \alpha(g) & \ast\\0 & \alpha(g)\end{pmatrix}$. Note that $\gamma$ sends the complex conjugation to $-1$, so its image has even order. 

Let $h \in K$. Then $h$ acts on $T_g/\mathfrak{p}$ with characteristic polynomial $(X-u)^2$, for some $u \in (\OO_{L_0}/\mathfrak{p})^{\times}$ with Teichm\"uller lift $u_0 \in \OO_{(L_0)_{\mathfrak{p}}}^{\times}$. Then $h' := u_0^{-1}h \in \GL{\OO_{(L_0)_{\mathfrak{p}}}}$ has finite order and its characteristic polynomial is congruent modulo $\mathfrak{p}$ to $(X-1)^2$. Therefore, the order of $h'$ is a power of $p$, so that $h^{p^n}$ is scalar for some $n \geq 1$. 

Hence the image of $K$ in $\PGL{\OO_{L_0}/\mathfrak{p}}$ is a $p$-group, and $\tilde{\rho}_g(G_{\Q})$ is an extension of the cyclic subgroup $C_1 := \gamma(G_{\Q})$ of order prime to $p$ by a $p$-group. Since $p > 5$, this is only possible if $\tilde{\rho}_g(G_{\Q})$ is a non-abelian dihedral group. In this situation, the surjection $\tilde{\rho}_g(G_{\Q}) \rar C_1$ factors through the group $\tilde{\rho}_g(G_{\Q})^{\mrm{ab}}$ of exponent $2$, so $C_1$ has order $2$. Therefore the kernel of $\tilde{\rho}_g(G_{\Q}) \rar C_1$ is a $p$-group with cardinality $\frac{Q}{2}$, so $\frac{Q}{2}$ is a power of $p$.

 }
 
\rem{Using the fact that $\SL{\OO_{F_{\mfk{p}}}}$ is its own derived subgroup, we can slightly refine the proof of Proposition \ref{only-we-se} to find $\tau \in G_{\Q}'$ acting on $T_{f,g,\mfk{p}}$ by $-\mrm{id}$. Let now $\eta: G_{\Q} \rar \{\pm 1\}$ be a quadratic character and assume that $T_{f,g,\mfk{p}}$ satisfies (wE) (resp. (sE)). Let $\sigma \in G_{\Q(\mu_{p^{\infty}})}$ be such that $\ker(\sigma-1 \mid T_{f,g,\mfk{p}}) \otimes (L_0)_{\mfk{p}}$ is a line (resp. $T_{f,g,\mfk{p}}/(\sigma-1)T_{f,g,\mfk{p}}$ is free of rank one over $\OO_{(L_0)_{\mfk{p}}}$). 

If $\eta(\sigma)=1$, then $\sigma$ also shows that $T_{f,g,\mfk{p}}\otimes \eta$ satisfies (wE) (resp. (sE)). If $\eta(\sigma)=-1$, then, since $\eta(\tau)=1$, $\sigma\tau$ shows that $T_{f,g,\mfk{p}} \otimes \eta$ satisfies (wE) (resp. (sE)). 

As a consequence, for any good prime $\mfk{p}$ for $(f,g)$, for any quadratic twists $f'$ of $f$ and $g'$ of $g$, $T_{f,g,\mfk{p}}$ satisfies (wE) (resp. (sE)) if and only if $T_{f',g',\mfk{p}}$ does. }

\medskip

We are now interested in the properties (wE) and (sE). Suppose that $\sigma \in G_{\Q(\mu_{p^{\infty}})}$ is such that $T_{f,g,\mfk{p}}/(\sigma-1)T_{f,g,\mfk{p}}$ is free of rank one over $\OO_{(L_0)_{\mfk{p}}}$. Let $M \in \GL{\OO_{(L_0)_{\mfk{p}}}}, N \in \GL{\OO_{L_0}}$ be the matrices of the actions of $\sigma$ on $T_{f,\mfk{p}}$ and $T_g$ respectively. Then $\det{M}$ and $N$ have finite order, and Corollary \ref{papier} shows that $M \in \alpha \GL{\OO_{F_{\mfk{p}}}}$ for a certain $\alpha$ contained in an explicit finite set of $\OO_{L,\mfk{p}}^{\times}$. It is thus natural to check which matrices $M \in \GL{\OO_{(L_0)_{\mfk{p}}}}, N \in \GL{\OO_{L_0}}$ are such that $N$ and $\det{M}$ have finite order and $\OO_{(L_0)_{\mfk{p}}}^{\oplus 4}/(M \otimes N-1)$ is a free $\OO_{(L_0)_{\mfk{p}}}$-module of rank one. 

\noindent
\textbf{Notation} \; Given $a,b \in R^{\times}$ for some ring $R$, let us write $\Delta_{a,b}=\begin{pmatrix}a & 0\\0 & b\end{pmatrix}$ and $U=\begin{pmatrix} 1 & 1\\0 & 1\end{pmatrix}$.  

The matrix $N$ has finite order, so it is similar to a diagonal matrix $\Delta_{u,v}$ where $u,v$ are roots of unity. In this case, after renaming if needed, it is necessary that $u \neq v$ and that $\ker(M-u^{-1}I_2)$ has rank one. Since $\det{M}$ is a fixed root of unity $\delta$, the two eigenvalues of $M$ are $u^{-1}$ and $\delta u$. The question then becomes: can we find a non-scalar $M \in \alpha \GL{\OO_{F_{\mfk{p}}}}$ such that its eigenvalues are $u^{-1}$ and $\delta u$? There is at least one case where we can answer in the affirmative: assuming that $u,\delta,\alpha \in F_{\mfk{p}}$, we can choose $M=\Delta_{u^{-1},\delta u}$ unless $u^{-1} =\delta u$, in which case we can take $M=u^{-1}U$.  

This is the rationale behind the following Proposition. 

\prop[sE-holds-in-splitting-set]{Assume that $\varepsilon_f\varepsilon_g \neq 1$. Then there exists a number field $F \subset L_1 \subset L_0$ such that, for all but finitely many good primes $\mathfrak{p} \subset \OO_{L_0}$, if the extension $L_1/F$ splits totally at the prime $\mathfrak{p} \cap \OO_F$, then $T_{f,g,\mathfrak{p}}$ has Euler-adapted Galois image. In particular, the set of such primes has positive lower density.}

\demo{By Proposition \ref{only-we-se}, we only need to show property (sE). 

The group $G_{\Q}$ is not covered by its proper subgroups $\ker{\varepsilon_f\varepsilon_g}$ and $\ker{\tilde{\rho}_g}$, so we can find $\sigma \in G_{\Q}$ such that $(\varepsilon_f\varepsilon_g)(\sigma) \neq 1$ and $\rho_g(\sigma)$ has distinct eigenvalues $u,v$. By the second part of Proposition \ref{jointimage-H}, we can choose $\sigma' \in G_{\Q(\mu_{p^{\infty}})} \cap \sigma G_H$ such that $\rho_g(\sigma')=\rho_g(\sigma)$, so, after replacing $\sigma$ with $\sigma'$, we may assume that $\sigma \in G_{\Q(\mu_{p^{\infty}})}$. Let $\alpha \in L^{\times}$ be attached to $\sigma$ by Corollary \ref{papier}. 

Let $L_1=F(\alpha,\varepsilon_f(\sigma),u)$. Let $p\nmid 30N_fN_g\varphi(N_f)\varphi(N_g)|\rho_g(G_{\Q})|$, where $\varphi$ denotes Euler's totient function, and let $\mathfrak{p}$ be a good prime with residue characteristic $p$ such that $L_1/F$ is totally split at $\mathfrak{p}$. We may asssume by Corollary \ref{papier} that $\alpha \in \OO_{L,\mfk{p}}^{\times}$, hence $\alpha \in \OO_{F_{\mfk{p}}}^{\times}$. 

Since $L_1/F$ is totally split at $\mfk{p}$, $\begin{pmatrix}\alpha & 0\\0 & \alpha^{-1}\varepsilon_f(\sigma)\end{pmatrix}\SL{\OO_{F_{\mfk{p}}}}$ is the set of matrices in $\GL{\OO_{F_{\mfk{p}}}}$ with determinant $\varepsilon_f(\sigma)$.

If $u^2 \neq \varepsilon_f(\sigma)^{-1}$, by Corollary \ref{papier}, we can choose $\sigma' \in \sigma G_{H(\mu_{p^{\infty}})}$ such that it acts on $T_{f,\mathfrak{p}}$ as $\Delta_{u^{-1},u\varepsilon_f(\sigma)}$ and on $T_g$ by $\rho_g(\sigma)$. Then the eigenvalues of $\sigma'$ on $T_{f,g,\mathfrak{p}}$ are $1,u^2\varepsilon_f(\sigma),(\varepsilon_f\varepsilon_g)(\sigma),v/u$: $1$ is distinct from all the others, hence (by our choice of $p$), not congruent to any of the others mod $\mathfrak{p}$, and $\sigma'$ works. 

If $u^2=\varepsilon_f(\sigma)^{-1}$, by Corollary \ref{papier}, we can choose $\sigma'\in \sigma G_{H(\mu_{p^{\infty}})}$ such that it acts on $T_{f,\mathfrak{p}}$ as $u^{-1}U$ and $\rho_g(\sigma)=\rho_g(\sigma')$: such a $\sigma'$ works. 
}
\medskip
\noindent

Proposition \ref{sE-holds-in-splitting-set} is not sufficient in order to answer Question \ref{loefflerqn}: in general, we cannot assume that we can find $\alpha,u$ that are both contained in $F_{\mfk{p}}$. There is however one special case in which the argument still applies: when $\alpha=1$ and $u = \pm 1$. The goal of the two following Propositions is to give sufficient conditions that imply that $\rho_g(G_H)$ contains a non-scalar matrix with $\pm 1$ as an eigenvalue.

\prop[special2]{%
Assume that $-1 \in \varepsilon_g(G_H)$ and that $\rho_g(G_H)$ does not contain an element similar to $\Delta_{-1,1}$.  
Then $g$ has either real or complex multiplication by some quadratic field $K$. If moreover we can choose $K$ so that it is not contained in $H$, then the following holds: 
\begin{itemize}[noitemsep,label=\tiny$\bullet$]
\item $i \in \varepsilon_g(G_H)$ and $\rho_g(G_H)$ contains $iI_2$,
\item the cardinality of $\tilde{\rho}_g(G_{\Q})$ is not divisible by $8$,
\item the cardinality of $\tilde{\rho}_g(G_{H})$ is not divisible by $4$.   
\end{itemize}
}

\demo{Let $M \in \GL{L_0}$ be a matrix with $\det{M}=-1$ and whose image in $\PGL{L_0}$ has order $2$. Then $M^2$ is scalar, so $M$ is conjugate in $\GL{\C}$ to $\Delta_{a,b}$ for some $a,b \in \C^{\times}$ with $a=-b$ and $ab=-1$. Therefore $M$ is similar to $\Delta_{1,-1}$ in $\GL{\C}$, so $M$ is similar to $\Delta_{1,-1}$ in $\GL{L_0}$.

Let $Z \leq \rho_g(G_H)$ be the subgroup of scalar matrices. Consider the surjective determinant map $\delta: \tilde{\rho}_g(G_H)=\rho_g(G_H)/Z \rar \varepsilon_g(G_H)/\det{Z}$, which has cyclic image. Let $x \in \tilde{\rho}_g(G_H)$ be an element of order $2$ such that $\delta(x)=-\det{Z}$: then there is a matrix $M \in \rho_g(G_H)$ above $x$ such that $\det{M}=-1$, so $M$ is similar to $\Delta_{1,-1}$, a contradiction. 

Clearly, $\rho_g(G_{\Q})'\subset \SL{L_0} \cap \rho_g(G_H)$, so $\delta$ vanishes on the subgroup $\tilde{\rho}_g(G_{\Q})'$ of $\tilde{\rho}_g(G_H)$. We distinguish two cases. \\

\emph{Case 1: when $-1 \notin \det{Z}$.} 

The image of $\delta$ contains exactly two elements of $2$-torsion: the trivial element $\det{Z}$ and $-\det{Z}$. For any $x \in \tilde{\rho}_g(G_H)$ of order two, $\delta(x) \in \mrm{im}(\delta)$ is a $2$-torsion element which is not $-\det{Z}$, so $x \in \ker{\delta}$. Since $\tilde{\rho}_g(G_H)/\ker{\delta}$ contains an element of order two, the subgroup of $\tilde{\rho}_g(G_H)$ generated by its elements of order two has even index: thus, $\tilde{\rho}_g(G_H)$ is not dihedral, and it is not isomorphic to $A_4$, $S_4$ or $A_5$ either. Since $\tilde{\rho}_g(G_H)$ is a subgroup of $\PGL{L_0}$, it is cyclic with order $n$ divisible by $4$. 

Therefore, the subgroup $\tilde{\rho}_g(G_{\Q})'$ of $\tilde{\rho}_g(G_H)$ is cyclic, so $\tilde{\rho}_g(G_{\Q})$ is not isomorphic to $A_4$, $S_4$, or $A_5$. Hence $\tilde{\rho}_g(G_{\Q})$ is dihedral and $g$ has either real or complex multiplication by a quadratic number field $K$. The dihedral group $\tilde{\rho}_g(G_{\Q})$ contains the cyclic subgroup $\tilde{\rho}_g(G_H)$ of order $n \geq 4$, so $\tilde{\rho}_g(G_H)$ is contained in the unique index two cyclic subgroup of $\tilde{\rho}_g(G_{\Q})$ (whose inverse image by $\tilde{\rho}_g$ is $G_K$), so that $K$ is unique and one has $K \subset H$, so we are done. 

\smallskip 

\emph{Case 2: when $-1 \in \det{Z}$.} 

In this case, $\rho_g(G_H)$ contains the scalar matrix $iI_2$. By the second paragraph of the proof, the group $\ker{\delta}$ does not contain an element of order $2$, so it has odd cardinality, and so does its subgroup $\tilde{\rho}_g(G_{\Q})'$. This is not possible if the group $\tilde{\rho}_g(G_{\Q})$ is isomorphic to $A_4$, $S_4$ or $A_5$ or if it is dihedral with cardinality divisible by $8$. Hence $\tilde{\rho}_g(G_{\Q})$ is dihedral with cardinality not divisible by $8$, so in particular, $g$ has either real or complex multiplication by some quadratic number field $K$. 

Let us furthermore assume that $K \not\subset H$. Let $C$ denote the cyclic subgroup of $D:= \tilde{\rho}_g(G_{\Q})$ such that $G_K=\tilde{\rho}^{-1}(C)$, and $D_H$ denote the subgroup $\tilde{\rho}_g(G_H)$ of $D$. Since $K \not\subset H$, one has $D_H \not\subset C$, so $D_H$ is dihedral with cyclic subgroup $D_H \cap C$ of index $2$ (it is possible, given our assumptions, that $D_H \simeq \Z/2\Z$ and $D_H \cap C$ is trivial; this case does not require any special treatment). In particular, the abelianization of $D_H$ has exponent $2$, so its cyclic quotient $\mrm{im}(\delta)$ has at most two elements. Hence either one has $\det{Z} = \varepsilon_g(G_H)$ or one has $\det{Z}=\varepsilon_g(G_H)^2$. 

Let $x \in D_H \backslash C$: then $x$ has order $2$, so $\delta(x) \neq -\det{Z}=\det{Z}$, so the image of $\delta$ has two elements, hence $-1 \in \det{Z} = \varepsilon_g(G_H)^2$, i.e. $i \in \varepsilon_g(G_H)$. Furthermore, $\ker{\delta}$ is a subgroup of index $2$ of $D_H$ which has odd cardinality, so the cardinality of $\tilde{\rho}_g(G_H)$ is not divisible by $4$. }

\smallskip

\medskip
\noindent

\prop[specialq]{Assume that there exists an odd prime $q \mid |\tilde{\rho}_g(G_{\Q})|$ such that $\varepsilon_g(G_H)$ contains an element of order $q$. Then one of the following statements holds:
\begin{itemize}[noitemsep,label=\tiny$\bullet$]
\item $q=3$ divides $[H:\Q]$, $\tilde{\rho}_g(G_{\Q})$ is isomorphic to $A_4$, $\tilde{\rho}_g(G_H)$ is isomorphic to $\F_2^{\oplus 2}$. 
\item $q=3$, $\tilde{\rho}_g(G_H) = \tilde{\rho}_g(G_{\Q})$ is isomorphic to $A_4$, $e^{2i\pi/9} \in \varepsilon_g(G_H)$, and the projection $\rho_g(G_H) \cap \SL{L_0} \rar \tilde{\rho}_g(G_H)$ is not surjective. 
\item $\rho_g(G_H)$ contains a matrix whose eigenvalues are $1$ and a primitive $q$-th root of unity.   
\end{itemize}}

\demo{As above, let $Z \leq \rho_g(G_H)$ be the subgroup of scalar matrices. 

Let us temporarily assume that $Z$ contains the $q$-th roots of unity, and suppose that there exists a matrix $M \in \rho_g(G_H) \cap \SL{L_0}$ whose image in $\PGL{L_0}$ has order $q$. Then $M^q = \{\pm I_2\}$, so that $M \sim t\Delta_{u,u^{-1}}$, where $t$ is a sign, and $u$ is a primitive $q$-th root of unity. Since $uI_2 \in Z$, $(uM)^2 \sim \Delta_{u^4,1}$ is also in $\rho_g(G_H)$.   

We now show that, except in specific situations which we can treat separately, $Z$ contains all the $q$-th roots of unity, and $\rho_g(G_H) \cap \SL{L_0}$ contains an element whose projective image has order $q$. 

As in the proof of Proposition \ref{special2}, the determinant map $\delta: \rho_g(G_H) \rar \varepsilon_g(G_H)/\det{Z}$ is surjective and its kernel contains $Z$ and $\rho_g(G_{\Q})'$. Hence $\delta$ factors through the quotient $\rho_g(G_H)/Z\rho_g(G_{\Q})' \simeq \tilde{\rho}_g(G_H)/\tilde{\rho}_g(G_{\Q})'$: in other words, $\varepsilon_g(G_H)/\det{Z}$ identifies as a cyclic quotient of $\tilde{\rho}_g(G_H)/\tilde{\rho}_g(G_{\Q})'$. \\

\emph{Case 1: $[\tilde{\rho}_g(G_{\Q}):\tilde{\rho}_g(G_{\Q})']$ is prime to $q$.} 

In this case, the cardinality of $\tilde{\rho}_g(G_{\Q})'$ is divisible by $q$, so $\tilde{\rho}_g(G_{\Q})'$ contains an element of order $q$. Now, $\tilde{\rho}_g(G_{\Q})'$ is the image in $\PGL{L_0}$ of $\rho_g(G_{\Q})' \subset \rho_g(G_H) \cap \SL{L_0}$ (the inclusion holds because $H/\Q$ is abelian). Hence the image of $\rho_g(G_H) \cap \SL{L_0}$ in $\PGL{L_0}$ contains an element of order $q$. 

Moreover, $\varepsilon_g(G_H)/\det{Z}$ is a quotient of the subgroup $\tilde{\rho}_g(G_{H})/\tilde{\rho}_g(G_{\Q})'$ of $\tilde{\rho}_g(G_{\Q})/\tilde{\rho}_g(G_{\Q})'$, so it has cardinality coprime to $q$. Since $\varepsilon_g(G_H)$ contains a primitive $q$-th root of unity $\omega$, one has $\omega \in \det{Z}$, so $Z$ contains a matrix $M$ with determinant $\omega$. Then $M^2=\omega I_2 \in Z$, and we are done. 

\smallskip

\emph{Case 2: $[\tilde{\rho}_g(G_{\Q}):\tilde{\rho}_g(G_{\Q})']$ is divisible by $q$.}

If $\tilde{\rho}_g(G_{\Q})$ is isomorphic to $S_4, A_5$ or is dihedral, then $\tilde{\rho}_g(G_{\Q})/\tilde{\rho}_g(G_{\Q})'$ is a $2$-group, so its cardinality is not divisible by $q$. Therefore $\tilde{\rho}_g(G_{\Q})$ is isomorphic to $A_4$, so $\tilde{\rho}_g(G_{\Q})/\tilde{\rho}_g(G_{\Q})'$ is isomorphic to $\Z/3\Z$, hence one has $q=3$.

Because $H/\Q$ is abelian, one has $\tilde{\rho}_g(G_{\Q})' \subset \tilde{\rho}_g(G_H) \subset \tilde{\rho}_g(G_{\Q})$, so either $\tilde{\rho}_g(G_H)=\tilde{\rho}_g(G_{\Q})$, or $\tilde{\rho}_g(G_H) =\tilde{\rho}_g(G_{\Q})' \simeq (\Z/2\Z)^{\oplus 2}$. By the argument at the beginning of the proof, it is enough to treat the following three cases. \smallskip

\emph{Case 2 (a): $\tilde{\rho}_g(G_H) = \tilde{\rho}_g(G_{\Q})$ and $Z$ does not contain any element of order $3$.}

Let $M \in \rho_g(G_H)$ be such that $\det{M}$ is a primitive third root of unity. If $M^2 \in Z$, $M^4=\det(M^2)I_2 \in Z$ is an element of order $3$, a contradiction: hence the projective image of $M$ has order at least $3$. The elements of $A_4$ have order $1$, $2$ or $3$, so the projective image of $M$ has order $3$, i.e. $M^3 \in \SL{L_0}$ is scalar, thus $M^3 \in \{\pm I_2\}$. Replacing $M$ with $M^2$ if needed, we may assume that $M^3=I_2$. Let $u,v \in L_0$ be the two eigenvalues of $M$, then one has $u^3=v^3=1$ and $u \notin \{v,v^{-1}\}$ (since $M$ is not scalar and $\det{M}$ is a primitive third root of unity), so exactly one of $u$ and $v$ is equal to $1$, and we are done. \smallskip

\emph{Case 2 (b): $\tilde{\rho}_g(G_H) = \tilde{\rho}_g(G_{\Q})$ and $e^{\frac{2i\pi}{9}} \notin \varepsilon_g(G_H)$.}

By Case 2 (a), we may assume that $Z$ (and therefore $\det{Z}$) contains an element of order $3$, so it is enough to show that there exists $M \in \rho_g(G_H) \cap \SL{L_0}$ with projective image of order $3$. The cyclic group $\varepsilon_g(G_H)$ contains no element of order $9$, so the quotient $\varepsilon_g(G_H)/\det{Z}$ has cardinality coprime to $3$. The abelianization of $\tilde{\rho}_g(G_H) \simeq A_4$ is isomorphic to $\Z/3\Z$, so $\delta$ is trivial and surjective, hence $\varepsilon_g(G_H)=\det{Z}$. This implies directly that $\rho_g(G_H) \cap \SL{L_0} \rar \tilde{\rho}_g(G_H)$ is surjective. We are done, because $\tilde{\rho}_g(G_H) \simeq A_4$ has an element of order $3$. \smallskip

\emph{Case 2 (c): $\tilde{\rho}_g(G_H)$ is a proper subgroup of $\tilde{\rho}_g(G_{\Q})$.}

In this case, $\tilde{\rho}_g(G_{\Q})/\tilde{\rho}_g(G_H) \simeq \Z/3\Z$ is a quotient of $\mrm{Gal}(H/\Q)$, so $3 \mid [H:\Q]$, and one has $\tilde{\rho}_g(G_H)=\tilde{\rho}_g(G_{\Q})' \simeq \F_2^{\oplus 2}$.  

}

\cor[sufficient-for-euler-adapted]{Let $\mfk{p}$ be a good prime for $(f,g)$. If $\rho_g(G_H)$ contains a non-scalar matrix $M$ such that $1$ or $-1$ is an eigenvalue of $M$, then $T_{f,g,\mathfrak{p}}$ satisfies (wE); if furthermore the order of $M$ is prime to $p$, then $T_{f,g,\mfk{p}}$ satisfies (sE). Hence $T_{f,g,\mfk{p}}$ has Euler-adapted Galois image for all but finitely many primes $\mathfrak{p}$. In particular, all of the above applies if any of the following conditions is satisfied:
\begin{enumerate}[noitemsep,label=(\roman*)]
\item \label{sc1} $\rho_g(G_H)$ contains the image under $\rho_g$ of the complex conjugation. This is the case if $H$ is unramified at every prime dividing the level of $g$. 
\item \label{sc2} $-1 \in \varepsilon_g(G_H)$ and $\tilde{\rho}_g(G_{\Q})$ is not dihedral.  
\item \label{sc3} $-1 \in \varepsilon_g(G_H)$, $\tilde{\rho}_g(G_H)$ is a dihedral group of order $4n$ with $n \geq 1$. 
\item \label{sc4} $-1 \in \varepsilon_g(G_H)$, $\tilde{\rho}_g(G_{\Q})$ is a dihedral group of order $8n$ with $n \geq 1$, and $g$ has either real or complex multiplication by some quadratic field $K \not\subset H$.
\item \label{sc5} $\varepsilon_g(G_H)$ contains $-1$ but not $i$, and, for some quadratic number field $K$ by which $g$ has either real or complex multiplication, one has $K \not\subset H$.
\item \label{sc6} The multiplicative order of $\varepsilon_f$ is a power of $2$, $\varepsilon_g(G_{\Q})$ contains an element of odd prime order $q \mid |\tilde{\rho}_g(G_{\Q})|$, and, if $\tilde{\rho}_g(G_{\Q}) \simeq A_4$, then $\rho_g(G_{\Q}) \cap \SL{L_0} \rar \tilde{\rho}_g(G_{\Q})$ is onto.  
\end{enumerate}}

\demo{By Proposition \ref{only-we-se}, we only need to discuss conditions (wE) and (sE).  
Let $\mfk{p}$ be a good prime with residue characteristic $p$ and $M \in \rho_g(G_H)$ be a non-scalar matrix admitting $s \in \{\pm 1\}$ as an eigenvalue. By Proposition \ref{jointimage-H}, there exists $\sigma \in G_{H(\mu_{p^{\infty}})}$ such that $\rho_g(\sigma)=M$ and $\sigma$ acts on $T_{f,\mfk{p}}$ by $sU$. Then the kernel of the action of $\sigma-1$ on $T_{f,g,\mfk{p}} \otimes (L_0)_{\mfk{p}}$ is a line and $T_{f,g,\mfk{p}}$ satisfies (wE). If the order of $M$ is prime to $p$, then $\OO_{(L_0)_{\mfk{p}}}^{\oplus 2}/(M-sI_2)\OO_{(L_0)_{\mfk{p}}}^{\oplus 2}$ is a free $\OO_{(L_0)_{\mfk{p}}}$-module of rank one, so $T_{f,g,\mfk{p}}/(\sigma-1)T_{f,g,\mfk{p}}$ is a free $\OO_{(L_0)_{\mfk{p}}}$-module of rank one, thus $T_{f,g,\mfk{p}}$ satisfies (sE). 

All that remains to do is show that in conditions \ref{sc1}-\ref{sc6}, there exists a non-scalar matrix $M \in \rho_g(G_H)$ which has $1$ or $-1$ as an eigenvalue. 

For \ref{sc1}, this is the case since $\rho_g$ maps the complex conjugation to a matrix similar to $\Delta_{1,-1}$. If we only assume that $H$ is unramified at every prime dividing the level of $g$, let $K_g/\Q$ be the finite extension such that $G_{K_g} = \ker{\rho_g}$. Since $\rho_g$ is unramified at every prime not dividing the level of $g$, there is no prime number $q$ that ramifies in both $K_g$ and $H$. Thus, $K_g \cap H$ is a number field unramified at every finite place of $\Q$, so it is $\Q$. Hence $G_H \rar \mrm{Gal}(K_g/\Q)$ is surjective and $\rho_g(G_H)=\rho_g(G_{\Q})$ contains the image of the complex conjugation.  

For \ref{sc2}, \ref{sc4}, \ref{sc5}, the conclusion follows from applying the contrapositive of Proposition \ref{special2}.   

Let us assume \ref{sc3} and let $K/\Q$ be a quadratic number field such that $g$ has either real or complex multiplication by $K$. Then $\tilde{\rho}_g(G_H)$ is not cyclic, so it is not contained in $\tilde{\rho}_g(G_K)$, so one has $K \not\subset H$. We can then apply the contrapositive of Proposition \ref{special2}.

For \ref{sc6}, let $r \geq 1$ be such that $\varepsilon_f^{2^r}$ is the trivial character. Let $(\gamma,\chi)$ be such that $f$ has an inner twist by $(\gamma,\chi)$. Then, by Lemma \ref{momose15}, $\chi$ is the product of a quadratic character and a power of $\varepsilon_f$, so $\chi^{2^r}$ is trivial, so $\mrm{Gal}(H/\Q)$ has exponent dividing $2^r$. Therefore, $d := [H:\Q]$ is a power of $2$. In particular, $\varepsilon_g(G_{\Q})/\varepsilon_g(G_H)$ and $\tilde{\rho}_g(G_{\Q})/\tilde{\rho}_g(G_H)$ are $2$-groups, so $\varepsilon_g(G_H)$ contains an element of order $q$ and $q \mid |\tilde{\rho}_g(G_H)|$.  

We want to apply Proposition \ref{specialq}, but we need to rule out its exceptional cases. The first one cannot occur since, under our assumptions, $[H:\Q]$ is not divisible by $3$. Now, assume that $\tilde{\rho}_g(G_{\Q})$ is isomorphic to $A_4$: in particular, it has no non-trivial quotient which is a $2$-group. The image $I$ of $\rho_g(G_H) \cap \SL{L_0} \rar \tilde{\rho}_g(G_{\Q})$ is a normal subgroup of $\tilde{\rho}_g(G_{\Q})$. Since $\rho_g(G_{\Q}) \cap \SL{L_0} \rar \tilde{\rho}_g(G_{\Q})$ is onto, $d$ kills the quotient $\tilde{\rho}_g(G_{\Q})/I$, so $\tilde{\rho}_g(G_{\Q})/I$ is a $2$-group. Hence it is trivial, i.e. $\rho_g(G_H) \cap \SL{L_0} \rar \tilde{\rho}_g(G_{\Q})$ is onto: the second exception of Proposition \ref{specialq} cannot occur and we are done. }

\rem{Because $G_H \subset \ker{\varepsilon_f}$, in all the cases of Corollary \ref{sufficient-for-euler-adapted}, the character $\varepsilon_f \varepsilon_g$ is not trivial. Corollary \ref{sufficient-for-euler-adapted} completes the proof of Theorem \ref{weakerconds}: its first two conditions follow from \ref{sc1}, the third condition from \ref{sc2}, the fourth condition from \ref{sc5}. The fifth condition follows in fact from the fourth: in this case, for every $(\gamma,\chi) \in \Gamma_f$, $\chi$ is quadratic by Lemma \ref{momose15}, and the forgetful map $\iota (\gamma,\chi) \in \Gamma_f \mapsto \chi \in \mrm{Hom}(G_{\Q},\{\pm 1\})$ is an injective group homomorphism. In particular, $\mrm{Gal}(H/\Q)$ is a $\F_2$-vector space, and, given a character $\psi: G_{\Q} \rar \{\pm 1\}$, $\psi \in \mrm{im}(\iota)$ if and only if $\psi$ vanishes on $\cap_{\chi\in \mrm{im}(\iota)}{\ker{\chi}}$.  }

\medskip

Let $E/\Q$ be an elliptic curve and $\rho: G_{\Q} \rar \GL{\C}$ be an odd Artin representation, with respective conductors $N_E,N_{\rho}$ and respective $L$-series $\sum_{n \geq 1}{a_n(E)n^{-s}}, \sum_{n \geq 1}{a_n(\rho)n^{-s}}$ (where $a_{\ell}(\rho)$ is the trace of the arithmetic Frobenius at $\ell$, for all but finitely many primes $\ell$), and if $\varepsilon$ denotes the primitive Dirichlet character attached to $\det{\rho}: G_{\Q} \rar \C^{\times}$, we define the \emph{Rankin--Selberg imprimitive $L$-series} $L(E,\rho,s)$ by 

\[L(E,\rho,s)=L^{(N_EN_{\rho})}(\varepsilon,2s-1)\sum_{n \geq 1}{a_n(E)a_n(\rho)n^{-s}},\] which is absolutely convergent for $\mrm{Re}(s) > \frac{3}{2}$ and defines a holomorphic function in this half-plane (where the exponent in the $L$-function of $\varepsilon$ mean that we remove the Euler factor at the places dividing $N_EN_{\rho}$). It is well-known (see \cite[\S 2.7]{KLZ15}) that this function extends to a meromorphic function over $\C$. 

\cor{Let $E$ be an elliptic curve over $\Q$ without complex multiplication and $\rho$ be a two-dimensional odd irreducible Artin representation of $G_{\Q}$, with splitting field $K$. Let $L_0/\Q$ be a finite extension containing all the eigenvalues of Frobenius for $\rho$ and $\mathfrak{p}$ be a prime ideal of $\OO_{L_0}$ with residue characteristic $p$. Let $\OO$ denote the valuation ring of $(L_0)_{\mfk{p}}$. Suppose that the following technical hypotheses hold:
\begin{enumerate}[noitemsep,label=(\roman*)]
\item $p$ is coprime to $30N_{\rho}N_E$, where $N_E$ (resp. $N_{\rho}$) is the conductor of $E$ (resp. $\rho$). 
\item If $\rho$ has projective dihedral image, then the cyclic subgroup is not a $p$-group.  
\item The $p$-adic Tate module of $E$ is a surjective representation $G_{\Q} \rar \GL{\Z_p}$. 
\item $E$ is ordinary at $p$ and $\rho(\Fr_p)$ has distinct eigenvalues mod $\mathfrak{p}$.
\end{enumerate} 
If $L(E,\rho,1) \neq 0$, then the group $\mrm{Hom}_{\OO[\mrm{Gal(K/\Q)}]}(\rho,\mrm{Sel}_{p^{\infty}}(E/K) \otimes_{\Z_p} \OO)$ is finite. }

\demo{This is \cite[Theorem 11.7.4]{KLZ15} without the assumption on the coprimality of conductors. Let $f \in \mathcal{S}_2(\Gamma_0(N_E))$ and $g \in \mathcal{S}_1(\Gamma_1(N_{\rho}))$ be the normalized newforms associated to $E$ and $\rho$: then $L(E,\rho,1)$ is the imprimitive Rankin--Selberg $L$-value $L(f,g,\mathbf{1},1)$ defined in \cite[\S 2.7]{KLZ15}. We can apply without any modifications the proof of \cite[Theorems 11.7.3, 11.7.4]{KLZ15}, as long as the following facts hold: 

\begin{enumerate}[noitemsep,label=(\alph*)]
\item $T_p E \otimes \rho \simeq T_{f,g,\mathfrak{p}}$ has Euler-adapted Galois image.
\item The $p$-adic $L$-value $L_p(f,g,1)$ (defined in \cite[Theorem 2.7.4]{KLZ15}) does not vanish. 
\end{enumerate} 

Since the modular form $f$ associated to the elliptic curve $E$ has no inner twists (since it has rational coefficients), one has $H=\Q$, so Proposition \ref{only-we-se} and Corollary \ref{sufficient-for-euler-adapted} show that $T_{f,g,\mathfrak{p}}$ has Euler-adapted Galois image. 

Now, we check the non-vanishing of $L_p(f,g,1)$. Let $\alpha_f,\beta_f \in (L_0)_{\mfk{p}}$ (resp. $\alpha_g,\beta_g \in (L_0)_{\mfk{p}}$) be the roots of the Hecke polynomial $X^2-a_p(f)X+p$ (resp. $X^2-a_p(g)X+\varepsilon(p)$), where $\alpha_f \in \OO_{(L_0)_{\mfk{p}}}^{\times}$.   
By the interpolation property of \cite[Theorem 2.7.4]{KLZ15}, $L_p(f,g,1)$ does not vanish as long as all of the following conditions are satisfied:
\[\beta_f \neq \alpha_f,\qquad \beta_f \neq p\alpha_f,\qquad \alpha_f \notin \{\alpha_g^{-1},\beta_g^{-1}\}, \qquad \beta_f \notin \{p\alpha_g^{-1},p\beta_g^{-1}\}.\]

Since $\alpha_f \in \OO_{(L_0)_{\mfk{p}}}^{\times}$ and $\alpha_f\beta_f=p$, one has $\alpha_f \neq \beta_f$. Furthermore, $\alpha_g, \beta_g$ are roots of unity, while, for any morphism $\sigma: (L_0)_{\mfk{p}} \rar \C$ of $L_0$-algebras, one has $|\sigma(\alpha_f)|=|\sigma(\beta_f)|=\sqrt{p}$, so the other conditions above are satisfied. }

\rem{When checking the Bloch--Kato conjecture, it seems more appropriate to use the $L$-function $L^{aut}(E,\rho,s)$ attached to the Rankin--Selberg convolution of the unitary automorphic representations attached to $f$ and $g$ (as in \cite[Section 3.6]{Bump}), since it admits a holomorphic continuation and satisfies a functional equation (with center $\frac{1}{2}$, see \cite[Theorem 11.7.1]{GH}). We are using the imprimitive $L$-function in the statement above primarily for the sake of simplicity, since it is the one appearing throughout \cite{KLZ15}. In fact, using the classical description of $L^{aut}(E,\rho,s)$ given in \cite[Theorem 2.2]{Li-RS}, it is not difficult to check that the meromorphic continuation of $L(E,\rho,s)$ is indeed holomorphic near $s=1$ and that one has $L(E,\rho,1) \neq 0$ if and only if $L^{aut}(E,\rho,\frac{1}{2}) \neq 0$. }

\section{Families of dihedral counter-examples: proof of Theorem \ref{inf-cex}}
\label{sect-theoD}

In this Section, we keep the notations from Section \ref{sect-theoA}, but furthermore assume that $\varepsilon_g=\varepsilon_f^{-1}\varepsilon$, where $\varepsilon$ is the character of some quadratic number field $K \subset H$ such that $g$ has either real or complex multiplication by $K$: in particular, $\varepsilon$ factors through $\mrm{Gal}(H/\Q)$. By construction, one has $G_H \subset \ker{\varepsilon_g}$, so the results of Propositions \ref{special2} and \ref{specialq} do not apply. One can nevertheless ask the following question: given a good prime ideal $\mathfrak{p}$ of $\OO_{L_0}$ for the couple $(f,g)$, when does $T_{f,g,\mathfrak{p}}$ satisfy (wE)? 

We first recall the following fact, already identified by Loeffler as an obstruction in \cite[Proposition 4.1.1]{bigimage}. 

\lem[efeg=1]{Let $M,N \in \GL{K}$ be matrices over a field $K$ such that $\det{M}\det{N}=1$. Then $\ker\left[M \otimes N-1\right]$ is never a line.}

\demo{This is clear if $M$ or $N$ is scalar, so we may assume that neither $M$ nor $N$ is scalar and that $1$ is an eigenvalue of $M \otimes N$. As in Section \ref{sect-theoA}, let $U$ denote the matrix $\begin{pmatrix}1 & 1\\0 & 1\end{pmatrix}$. 

Assume that $M \sim \alpha U$ for some $\alpha \in K^{\times}$: then $\alpha^{-1}$ is an eigenvalue of $N$ and $\det{M}=\alpha^2$. Since $\det{M}\det{N}=1$, the characteristic polynomial of $N$ is $(X-\alpha^{-1})^2$. The matrix $N$ is not scalar, so one has $N \sim \alpha^{-1} U$. It is easy to see that $U \otimes U - 1$ has rank two, so $\ker(M \otimes N-1)$ is not a line. 

Thus, we can (and do) assume that $M$ and $N$ both have distinct eigenvalues, respectively $m,m'$ and $n,n'$. Then $M \otimes N$ is diagonalizable with eigenvalues $mn,m'n',m'n,mn'$. Since $1=(mn)(m'n')=(mn')(m'n)$, there is never exactly one of these eigenvalues which is equal to one.}

\medskip

As noted by Loeffler, this lemma implies that, for any couple $(f',g')$ of newforms, $T_{f',g',\mathfrak{p}}$ cannot satisfy (wE) if the product of the characters of $f',g'$ is trivial. That is why it is assumed in Question \ref{loefflerqn} that $\varepsilon_f\varepsilon_g \neq 1$. 

\rems{
\begin{enumerate}[noitemsep,label=(\roman*)]
\item \label{starts3-1} If $\eta$ is a quadratic Dirichlet character, then the newform $f \otimes \eta$ has the same weight and character as $f$ and the same attached fields $H,L,F$. Hence, if $\beta$ is another quadratic Dirichlet character, the couple of newforms $(f \otimes \eta, g \otimes \beta)$ also satisfies the assumptions made at the beginning of this Section. 
\item \label{starts3-2} Let $\mfk{p} \subset \OO_{L_0}$ be a prime ideal and $\sigma \in G_{\Q}$ be such that the kernel of $\sigma-1$ acting on $T_{f,g,\mfk{p}}$ is a line. Let $M \in \GL{\OO_{(L_0)_{\mfk{p}}}}$ and $N \in \GL{\OO_{L_0}}$ be the matrices of the action of $\sigma$ on $T_{f,\mfk{p}}, T_{g}$. Then one has $\varepsilon(\sigma)=\varepsilon_f(\sigma)\varepsilon_g(\sigma) = \det{M}\det{N}\neq 1$ by Lemma \ref{efeg=1}, so $\sigma \notin G_K$. Since $g$ has real or complex multiplication by $K$, $N$ has trace zero, so, if $u \in \C^{\times}$ denotes a square root of $-\det{N}$, $N$ is similar to $\Delta_{u,-u}$. This implies that $N^2$ is a scalar matrix. Moreover, $N$ is similar to $-N$, so $M \otimes N$ is similar to $M \otimes (-N)=-(M \otimes N)$, so that $\ker(M \otimes N+1)$ is also a line. 
\item \label{starts3-3} The discussion of \ref{starts3-2} has the following implication to the situation of \ref{starts3-1}: the representation $T_{f,g,\mfk{p}}$ satisfies (wE) if and only if $T_{f \otimes \eta, g \otimes \beta, \mfk{p}}$ satisfies (wE).
\end{enumerate}} 

\defi{Let $\mathcal{B}: \mrm{Gal}(H/\Q) \times G_F \rar \{\pm 1\}$ be the following map. Let $(\sigma,\tau) \in \mrm{Gal}(H/\Q) \times G_F$. By Corollary \ref{papier}, there is some $\alpha \in L^{\times}$ such that $\gamma(\alpha)=\chi(\sigma)\alpha$ for all $(\gamma,\chi) \in \Gamma_f$. Since $\gamma(\varepsilon_f)=\chi^2\varepsilon_f$, one has $\nu := \frac{\alpha^2}{\varepsilon_f(\sigma)} \in F^{\times}$. Then $\mathcal{B}(\sigma,\tau)=\frac{\tau(\sqrt{\nu})}{\sqrt{\nu}}$. }

\prop{$\mathcal{B}$ is a bilinear pairing, and factors through a finite quotient of $G_F$. Moreover, for every $\sigma \in \mrm{Gal}(H/\Q)$, $\mathcal{B}(\sigma, \cdot): G_F \rar \{\pm 1\}$ is unramified away from $2 \cdot \mrm{disc}(L/F)$.}

\demo{The bilinearity is a direct verification. Moreover, the subgroup $\varepsilon_f(G_{\Q})$ of $L^{\times}$ is finite cyclic, so it is generated by some root of unity $\zeta$. It is then easy to check that if $\tau \in \mrm{Gal}(\overline{F}/L(\sqrt{\zeta}))$, then one has $\mathcal{B}(\cdot,\tau)=1$.}

\medskip
\noindent
Let $\mathfrak{p} \subset \OO_{L_0}$ be a good prime ideal for $(f,g)$, we write $\mathcal{B}_{\mathfrak{p}}=\mathcal{B}(\cdot,\Fr_{\mathfrak{p}}) : \mrm{Gal}(H/\Q) \rar \{\pm 1\}$. Denote by $M$ the subgroup of $\mrm{Gal}(H/\Q)$ made by those $\sigma \in \mrm{Gal}(H/\Q)$ such that $\mathcal{B}(\sigma,\cdot)=1$. By Chebotarev, $M$ is thus the intersection of the kernels of the $\mathcal{B}_{\mathfrak{p}}$.

The motivation for the construction of $\mathcal{B}$ is the following.

\prop[se-equiv-Bp]{Let $\mathfrak{p}$ be a prime ideal of $\OO_{L_0}$ of residue characteristic $p$, good for $(f,g)$. Then the following are equivalent: 
\begin{itemize}[noitemsep, label=\tiny$\bullet$]
\item $T_{f,g,\mathfrak{p}}$ satisfies (wE),
\item $T_{f,g,\mathfrak{p}}$ satisfies (sE),
\item $\mathcal{B}_{\mathfrak{p}} \neq \varepsilon$. 
\end{itemize}}

\demo{Suppose that $T_{f,g,\mathfrak{p}}$ satisfies (wE). Let $\sigma \in G_{\Q(\mu_{p^{\infty}})}$ be such that $\ker(\sigma - 1 \mid T_{f,g,\mathfrak{p}})$ is a line. Let $\alpha \in \OO_{L,\mathfrak{p}}^{\times}$ be associated to $\sigma$ by Corollary \ref{papier}. Let $M \in \GL{\OO_{(L_0)_{\mathfrak{p}}}}, N \in \GL{\OO_{L_0}}$ be the matrices given by the action of $\sigma$ on $T_{f,\mathfrak{p}}$ and $T_g$, respectively. By point \ref{starts3-2} in the Remarks following Lemma \ref{efeg=1}, one has $\sigma \notin G_K$ and the matrix $N^2$ is scalar. In fact, one has by Cayley--Hamilton $N^2=(-\det{N})I_2 = (-\varepsilon_g(\sigma))I_2 = (-\varepsilon(\sigma)\varepsilon_f(\sigma)^{-1})I_2=\varepsilon_f(\sigma)^{-1}I_2$. Let $\omega \in L_0$ be a square root of $\varepsilon_f(\sigma)^{-1}$, then the two eigenvalues of $N$ are $\omega$ and $-\omega$. 

By Corollary \ref{papier}, we can write $M=\alpha^{-1}\varepsilon_f(\sigma)M'$, where $M' \in \GL{\OO_{F_{\mathfrak{p}}}}$ has determinant $\frac{\alpha^2}{\varepsilon_f(\sigma)}$. Thus $\alpha\omega^2$ is an eigenvalue of $M' \otimes \Delta_{\omega,-\omega}$, so that one of the eigenvalues of $M'$ is $\pm \omega \alpha$. Since $(\pm \omega \alpha)^2=\det{M'}$, one has $\pm \alpha \omega = \frac{1}{2}\Tr{M'} \in F_{\mathfrak{p}}$, so $\mathcal{B}(\sigma,\Fr_{\mathfrak{p}})=1 \neq \varepsilon(\sigma)$.

Suppose conversely that $\mathcal{B}(\cdot,\Fr_{\mathfrak{p}})\neq \varepsilon$. Since $\mathcal{B}_{\mfk{p}}^2=1$ and $\varepsilon \neq 1$, $\varepsilon$ is not a power of $\mathcal{B}_{\mfk{p}}$, so $\ker{\mathcal{B}(\cdot,\Fr_{\mathfrak{p}})}$ is not contained in $\ker{\varepsilon}$. Hence, there exists $s \in G_{\Q}$ such that $\varepsilon(s)=-1$ and $\mathcal{B}_{\mathfrak{p}}(s)=1$. Let $\alpha \in \OO_{L,\mathfrak{p}}^{\times}$ be attached to $s$ by Corollary \ref{papier} and $\nu := \frac{\alpha^2}{\varepsilon_f(s)} \in \OO_{F,\mfk{p}}^{\times}$. Since $\mathcal{B}_{\mathfrak{p}}(s)=1$, $\nu$ is the square of some $\nu' \in \OO_{F_{\mathfrak{p}}}^{\times}$. Therefore, $\varepsilon_f(s)^{-1}=\nu\alpha^{-2}$ is the square of $\omega := \nu'\alpha^{-1} \in L_{\mathfrak{p}}$.  

By Proposition \ref{jointimage-H} and Corollary \ref{papier}, we can choose $\sigma_1 \in G_{\Q(\mu_{p^{\infty}})} \cap sG_H$ acting on $T_{f,\mathfrak{p}}$ by the matrix $\alpha^{-1}\varepsilon_f(s)\begin{pmatrix}\nu' & 1\\0 & \nu'\end{pmatrix}$. Since $\varepsilon(\sigma_1)=\varepsilon(s)=-1$, the matrix $\tilde{\rho}_g(\sigma_1) \in \PGL{\OO_{L_0}}$ has order two. Since $\det{\rho_g(\sigma_1)}=-\varepsilon_f(\sigma_1)^{-1}=-\varepsilon_f(s)^{-1}$ and $p \neq 2$, $\rho_g(\sigma_1)$ is $\GL{\OO_{(L_0){\mathfrak{p}}}}$-conjugate to $\Delta_{\omega,-\omega}$. Therefore, $\sigma_1$ acts on $T_{f,g,\mathfrak{p}}$ by $\begin{pmatrix}1 & a & 0 & 0\\0 & 1 & 0 & 0\\0 & 0 & -1 & \ast\\0 & 0 & 0 & -1\end{pmatrix}$ with $a \in \OO_{(L_0)_{\mathfrak{p}}}^{\times}$: hence $T_{f,g,\mathfrak{p}}/(\sigma_1-1)T_{f,g,\mathfrak{p}}$ is free of $\OO_{(L_0)_{\mathfrak{p}}}$-rank one. 
}

\cor[bidual]{The answer to Question \ref{loefflerqn} for $(f,g)$ is negative if, and only if, $\varepsilon(M)=1$.}

\demo{By Chebotarev and Proposition \ref{se-equiv-Bp}, the answer to Question \ref{loefflerqn} for $(f,g)$ is positive if and only if $\varepsilon$ cannot be written as a $\mathcal{B}(\cdot,\tau)$ for any $\tau \in G_F$. Because $\mathcal{B}$ is bilinear with values in a $\F_2$-vector space, the set $\{\mathcal{B}(\cdot,\tau) \mid \tau \in G_F\}$ is exactly the set of homomorphisms $f: \mrm{Gal}(H/\Q) \rar \{\pm 1\}$ such that $f(M)=1$, and the conclusion follows.}

\rem{It is a direct verification that $\mathcal{B}$ does not change if $f$ is replaced with one of its quadratic twists. Thus, at least when $\mfk{p}$ has large enough residue characteristic, the equivalence of Proposition \ref{se-equiv-Bp} does not change if $f$ is replaced with one of its quadratic twists. }

\prop[whenepsilonexists]{Let $H_0 \subset H$ be the subfield such that $G_{H_0} = \ker{\varepsilon_f}$. Let $c \in \mrm{Gal}(H/\Q)$ denote the complex conjugation. 
Assume that one of the following holds:
\begin{enumerate}[noitemsep,label=$(\alph*)$]
\item\label{weps-a} $k$ is even and $H$ is not totally real, 
\item\label{weps-b} $k$ is odd and $[H:H_0]>2$,  
\item\label{weps-c} $k$ is odd, $[H:H_0]=2$, and $\varepsilon_f$ does not take its values in $L^{\times 2}$ (this last condition is satisfied if $L$ does not contain a fourth root of unity $i$), 
\item\label{weps-d} $k$ is odd, $[H:H_0]=2$, $i \in L$, and for all $(\gamma,\delta) \in \Gamma_f$, $\gamma(i)=\delta(-1)i$,
\item\label{weps-e} $k$ is odd, $\mrm{Gal}(H/H_0)$ contains a unique non-trivial element $\sigma_0$, $i \in L$ and, for some $(\gamma,\delta) \in \Gamma_f$, one has $\gamma(i) \neq \delta(\sigma_0)\delta(-1)i$.
\end{enumerate}
Then there exists a non-trivial character $\varepsilon: \mrm{Gal}(H/\Q)/M \rar \{\pm 1\}$ sending $c$ to $(-1)^{k-1}$. \\
Moreover, when $k$ is odd, $[H:H_0]=2$, $i \in L$ but is not a value of $\varepsilon_f$, then the existence of such a $\varepsilon$ is equivalent to the last clause of \ref{weps-e}. Under the same assumptions, if $\varepsilon_f^2=1$, then \ref{weps-d} and \ref{weps-e} are equivalent.
}
 
\demo{Remember that $L$ is the number field generated by the coefficients of $f$, so it is a subfield of $\C$. Note that $\mathcal{B}$ realizes an isomorphism from $\mrm{Gal}(H/\Q)/M$ to a subgroup $B$ of the $\F_2$-vector space $\mrm{Hom}(G_F,\F_2)$, so $\mrm{Gal}(H/\Q)/M$ has exponent dividing $2$. Let $c \in \mrm{Gal}(H/\Q)$ be the image of the complex conjugation. 

Let $\sigma \in M$ be such that $\varepsilon_f(\sigma)=1$: let $\alpha \in L^{\times}$ be associated to $\sigma$ by Corollary \ref{papier} and $\nu=\frac{\alpha^2}{\varepsilon_f(\sigma)}=\alpha^2 \in F^{\times}$. Then, for all $\tau \in G_F$, $1=\mathcal{B}(\sigma,\tau)=\frac{\tau(\sqrt{\nu})}{\sqrt{\nu}}=\frac{\tau(\alpha)}{\alpha}$, hence $\alpha \in F^{\times}$. Therefore, for all $(\gamma,\delta) \in \Gamma_f$, $\alpha=\gamma(\alpha)=\delta(\sigma)\alpha$, so $\delta(\sigma)=1$, hence $\sigma \in G_H$. Thus, $(\mathcal{B},\varepsilon_f): \mrm{Gal}(H/\Q) \rar \mrm{Hom}(G_F,\F_2) \times \C^{\times}$ is injective. 

In particular, $\sigma \in \mrm{Gal}(H/H_0) \longmapsto \mathcal{B}(\sigma,\cdot)$ is injective. In particular, if $k$ is even and $H$ is not totally real, then $c$ is non-trivial: since $\varepsilon_f(c)=1$, $c \notin M$, so there is a morphism $\mrm{Gal}(H/\Q)/M \rar \{\pm 1\}$ sending $c$ to $-1$. 

Suppose from now on that $k$ is odd. The goal is to show that there exists a non-trivial character $\mrm{Gal}(H/\Q)/\langle M,c\rangle \rar \{\pm 1\}$. To do that, we need to show that $\mrm{Gal}(H/\Q) \supsetneq \langle M,c\rangle$, or in other words that $B \supsetneq \{1,\mathcal{B}(c,\cdot)\}$. Since $\sigma \in \mrm{Gal}(H/H_0) \longmapsto \mathcal{B}(\sigma,\cdot)$ is injective, one has $|B| \geq [H:H_0]$, so we are done if $[H:H_0] > 2$. 

Suppose from now on that $[H:H_0]=2$ and let $\sigma_0$ be the unique non-trivial element of $\mrm{Gal}(H/H_0)$. For any $\sigma \in \mrm{Gal}(H/\Q)$, $\sigma \in M$ if and only if the element $\alpha \in L^{\times}$ associated to $\sigma$ by Corollary \ref{papier} is such that $\frac{\alpha^2}{\varepsilon_f(\sigma)} \in F^{\times 2}$, which is the case if and only if $\sqrt{\varepsilon_f(\sigma)}$ lies in $L^{\times}$ and is attached to $\sigma$ by Corollary \ref{papier}.  

We have an exact sequence of Abelian groups \[0 \rar \mrm{Gal}(H/H_0) \rar \mrm{Gal}(H/\Q) \rar \mrm{Gal}(H_0/\Q) \rar 0.\] Since $[H:H_0]=2$, the group $\mrm{Gal}(H/H_0)$ is cyclic of order two. The group $\mrm{Gal}(H_0/\Q)$ is isomorphic to the image of $\varepsilon_f$ (which contains $-1$), so it is a cyclic group with even cardinality $n$. Fix $\sigma \in \mrm{Gal}(H/\Q)$ such that $\varepsilon_f(\sigma)$ has order $n$. Then $\sigma^n \in \mrm{Gal}(H/H_0) \cap M$ because $n$ is even and $\mrm{Gal}(H/\Q)/M$ has exponent dividing $2$, so $\sigma^n$ is trivial by the second paragraph. Thus $\sigma_{|H_0} \mapsto \sigma$ defines a splitting of the above exact sequence, and therefore $\mrm{Gal}(H/\Q) = \Z/2\Z\sigma_0 \oplus \Z/n\Z \sigma$. Hence, representatives for $\mrm{Gal}(H/\Q)/M$ are given (possibly with repetitions) by $\sigma_0,\sigma_0\sigma,\sigma,\mrm{id}$. If $\sigma$ (or $\sigma\sigma_0$) is in $M$, the previous paragraph shows that $\varepsilon_f(\sigma)=\varepsilon_f(\sigma\sigma_0) \in L^{\times 2}$. Since $\varepsilon_f(\sigma)$ is a generator of the image of $\varepsilon_f$, this implies that $\varepsilon_f$ takes its values in $L^{\times 2}$. Hence, in the case \ref{weps-c}, $M$ does not contain $\sigma$ or $\sigma_0\sigma$, so $|B| > 2$ and we are done. We assume from now on that $i \in L$. 

Note that \ref{weps-d} (resp. \ref{weps-e}) holds if and only if $c \in M$ (resp. $\sigma_0c \notin M$). Hence, since $\sigma_0 \notin M$, \ref{weps-d} implies \ref{weps-e}. Now, when \ref{weps-e} holds, one has $\sigma_0 c \notin M$ and $\sigma_0 \notin M$, so $\sigma_0$ is not in the subgroup generated by $c$ and $M$, and the conclusion follows. 

If $i$ is not a value taken by $\varepsilon_f$ (which is for instance the case if $\varepsilon_f^2=1$), then $n$ is not divisible by $4$. Because $\mrm{Gal}(H/\Q)$ is isomorphic to $\Z/2\Z \oplus \Z/n\Z$, the $2$-torsion subgroup of $\mrm{Gal}(H/\Q)$ is isomorphic to $\F_2^{\oplus 2}$ and projects isomorphically on $\mrm{Gal}(H/\Q) \otimes \F_2$, hence it is onto $\mrm{Gal}(H/\Q)/M$. Because $\varepsilon_f(c)=(-1)^k=-1$, $c,\sigma_0$ are two distinct non-trivial $2$-torsion elements in $\mrm{Gal}(H/\Q)$, so $(u,v) \in \mathbb{F}_2^{\oplus 2} \mapsto \sigma_0^uc^v \in \mrm{Gal}(H/\Q)/M$ is onto. Since $\sigma_0 \notin M$, $\mrm{Gal}(H/\Q)$ is then generated by $M$ and $c$ if and only if $\sigma_0 c \in M$, if and only if \ref{weps-e} does not hold, and the conclusion follows.

Suppose finally that $\varepsilon_f^2=1$. By Lemma \ref{momose15}, for any $(\gamma,\delta) \in \Gamma_f$, $\delta$ is a quadratic character, so $\mrm{Gal}(H/\Q)$ is a $\F_2$-vector space; by the previous paragraph, it has dimension two and $(\sigma_0,c)$ is a basis of this vector space. One can check directly that $\Gamma_f$ also has exponent two and that the bilinear pairing $[(\gamma,\delta),\sigma] \in \Gamma_f \times \mrm{Gal}(H/\Q) \mapsto \delta(\sigma) \in \{\pm 1\}$ is a perfect duality of $\F_2$-vector spaces. Therefore, the homomorphism $(\gamma,\delta) \in \Gamma_f \mapsto \frac{\gamma(i)}{i} \in \{\pm 1\}$ is given by $(\gamma,\delta) \mapsto \delta(\sigma_1)$ for some $\sigma_1 \in \{\sigma_0,c,\sigma_0c\}$. By taking $(\gamma,\delta)=(c,\varepsilon_f^{-1})$, we find that $\varepsilon_f(\sigma_1)=-1$, so $\sigma_1 \in \{c,\sigma_0c\}$. Then $i$ is the element $\alpha$ attached to $\sigma_1$ by Corollary \ref{papier}, so it follows that $\mathcal{B}(\sigma_1,\cdot)=1$ i.e. $\sigma_1 \in M$. Because $\sigma_0 \notin M$, $\sigma_1$ is the only element of $\{c,\sigma_0c\}$ which lies in $M$. Now, by definition, \ref{weps-d} holds if and only if $\sigma_1=c$, if and only if $\sigma_0 c \notin M$, if and only if \ref{weps-e} holds.

}

By Corollary \ref{bidual}, Proposition \ref{whenepsilonexists} implies: 

\cor{Suppose that the hypotheses of Proposition \ref{whenepsilonexists} hold: by Proposition \ref{whenepsilonexists}, there exists a non-trivial character $\varepsilon: \mrm{Gal}(H/\Q)/M \rar \{\pm 1\}$ mapping the complex conjugation to $(-1)^{k-1}$. Let $g$ be a weight one newform with character $\varepsilon\varepsilon_f^{-1}$ which has either complex multiplication by $\varepsilon$ if $k$ is even, or real multiplication by $\varepsilon$ is $k$ is odd. Then there are infinitely many maximal ideals $\mfk{p}\subset \mathcal{O}_{L_0}$ (with residue characteristic $p$) such that there is no $\sigma \in G_{\Q(\mu_{p^{\infty}})}$ such that the kernel of the action of $\sigma-1$ on $V_{f,g,\mfk{p}}$ is a line.

In other words, the answer to Question \ref{loefflerqn} for the couple $(f,g)$ of modular forms is negative.}

\rems{
\begin{enumerate}[noitemsep,label=(\roman*)]
\item Let $f \in \mathcal{S}_2(\Gamma_0(63))$ be a newform with LMFDB \cite{lmfdb} label $63.2.a.b$\footnote{Technically, these labels refer to a Galois orbit of newforms; however, the choice of embedding does not matter in all the situations below, so we omit it.}: it has coefficients in $\Q(\sqrt{3})$, an inner twist by the quadratic character $\varepsilon$ of $K=\Q(\sqrt{-3})$ and no complex multiplication. The newform $g \in \mathcal{S}_1(\Gamma_1(675))$ with LMFDB label $675.1.c.b$ has rational coefficients, character $\varepsilon$ and complex multiplication by $K$. The argument shows that (after possibly removing a finite number of primes) $T_{f,g,\mathfrak{p}}$ does not satisfy (wE) when the residue characteristic is congruent to $5,7$ mod $12$.
\item Consider a newform $f \in \mathcal{S}_3(\Gamma_1(75))$ with LMFDB label $75.3.d.b$. Then $L=\Q(i,\sqrt{5})$, $F=\Q$, $H_0=\Q(\sqrt{-15})$, $H=\Q(\sqrt{-3},\sqrt{5})$. Then it can be checked that $\mrm{Gal}(H/\Q)$ is generated by $M$ and $c$. When $\varepsilon$ is the character of the unique real quadratic field contained in $H$ and $g$ is a weight one newform with character $\varepsilon\varepsilon_f^{-1}$ and with real multiplication by $\varepsilon$, property (sE) then holds for all but finitely many $T_{f,g,\mathfrak{p}}$, even though Theorem \ref{weakerconds} does not apply.
\item When $k$ is even, $H_0$ is necessarily totally real. However, $H$ can be totally real and larger than $H_0$. This is what happens for the newform $f \in \mathcal{S}_2(\Gamma_0(289))$ with LMFDB label $289.2.a.f$: indeed, $H_0=\Q$, but $H=\Q(\sqrt{17})$. By Theorem \ref{weakerconds}, the answer to Question \ref{loefflerqn} is positive for $f$ and any weight one newform. 
\item For a newform $f \in \mathcal{S}_3(\Gamma_1(24))$ with LMFDB label $24.3.h.c$, the condition \ref{weps-c} of Proposition \ref{whenepsilonexists} is satisfied: indeed, in this situation $L=\Q(\sqrt{2},\sqrt{-7})$ does not contain $i$, while $H=\Q(\sqrt{2},\sqrt{-3})$ and $H_0=\Q(\sqrt{-6})$ satisfy $[H:H_0]=2$. Indeed, $B$ contains the two distinct non-trivial elements $\mathcal{B}(c,\cdot)$ and $\mathcal{B}(\Fr_5,\cdot)$ (they are the Kummer characters associated to $\sqrt{14}$ and $\sqrt{2}$ respectively), so $M$ is trivial and $\mrm{Gal}(H/\Q)$ properly contains $\{\mrm{id},c\}$.
\item For a newform $f \in \mathcal{S}_3(\Gamma_1(64))$ with LMFDB label $64.3.d.a$, we have $L=\Q(\mu_{12})$ and $F=\Q$; moreover, $H=\Q(\mu_8)$ and $\mrm{Gal}(H/H_0)$ corresponds to the subgroup $\{1,3\}$ of $(\Z/8\Z)^{\times}$. Because $a_{63}(f) \in iF^{\times}$, the condition \ref{weps-d} of Proposition \ref{whenepsilonexists} is verified.   
\item Even when $i \in L$ is not a value of $\varepsilon_f$, it is not true that \ref{weps-e} implies \ref{weps-d}. A counter-example is given by the newform $f \in \mathcal{S}_3(\Gamma_1(324))$ with LMFDB label $324.3.f.s$\footnote{With the tools dicussed here, this verification is not difficult using the LMFDB data and MAGMA \cite{magma}; we omit it for the sake of brevity.}. In this situation, $\varepsilon_f$ has conductor $36$ and order $6$, while $L$ has degree $16$ and $|\Gamma_f|=8$, so that $F=\Q(\sqrt{57})$. 
\item Consider a newform $f \in \mathcal{S}_3(\Gamma_1(21))$ with LMFDB label $21.3.h.b$. In this case, one has $|\Gamma_f|=4$, but $H=H_0=\Q(\sqrt{-3},\mu_7)^+$. In fact, $F=\Q$, $L=\Q(\sqrt{-3},\sqrt{-5})$, and one can directly compute that $B = \{\mathcal{B}(\mrm{id},\cdot),\mathcal{B}(c,\cdot)\}$.  
\end{enumerate}
}

\medskip
\noindent
The following result completes the proof of Theorem \ref{inf-cex}. 

\prop{Let $\varepsilon$ be the character of the imaginary (resp. real) quadratic field $K$, and $\chi$ be an even (resp. odd) primitive Dirichlet character. There exists an infinite collection $\mathcal{G}$ of weight one newforms $g$ with complex (resp. real) multiplication by $K$ and character $\varepsilon\chi$ such that no two elements of $\mathcal{G}$ are twists of each other by a quadratic character.  }

\demo{This statement is easier to establish using Galois representations. \medskip

\emph{Step 1: $\chi$ factors through the transfer map $V: G_{\Q}^{ab} \rar G_K^{ab}$.}

First, we show that there is a finite abelian extension $K'/K$ such that the kernel of the transfer map $V_{K'}: \mrm{Gal}(K'/\Q)^{ab} \rar \mrm{Gal}(K'/K)^{ab}$ is contained in that of $\chi$. 

Let $C$ be the conductor of $\chi$, $D$ be a positive multiple of $C$ and $K'/K$ the finite abelian extension whose group of norms is $K^{\times}\prod_v{M_v}$ where $M_v=\C^{\times}$ if $v$ is complex, $\R^{\times}_+$ if $v$ is real, $\OO_{K_v}^{\times}$ if $v$ is finite and not dividing $D$, and $1+D\OO_{K_v}$ otherwise. Because this group of norms is stable under $\mrm{Gal}(K/\Q)$, $K'/\Q$ is Galois. By functoriality in class field theory \cite[Proposition IV.5.9, Sections VI.4, VI.5]{NeuANT}, $\chi$ vanishes on the kernel of $V_{K'}$ as long as the following claim is true: for any $z \in \mathbb{A}_{\Q}^{\times}$ whose image in $\mathbb{A}_K^{\times}$ is contained in $K^{\times}\prod_v{M_v}$, we can write $z=\kappa(z_{\infty},z')$, where $\kappa \in \Q^{\times}$, $(z_{\infty}, z') \in \R^{\times} \times (\hat{\Z}^{\times} \cap (1+C\hat{\Z}))$, and $z_{\infty} > 0$ if $K$ is real. 

When $K$ is imaginary, we may assume (after multiplying $z$ with a rational number) that $z=(z_0,z_f)$ with $z_f \in \hat{\Z}^{\times}$ and $z_0 \in \R^{\times}_+$, and thus that its image in $\mathbb{A}_K^{\times}$ is contained in $\OO_K^{\times}\prod_{v}{M_v}$. It is then enough to find $D$ such that any unit $u \in \OO_K^{\times}$ which is congruent to an integer modulo $D$ is congruent to $\pm 1$ modulo $C$ (this $\pm 1$ corresponds to the fact that the complex conjugation is contained in the kernel). Since $\OO_K^{\times}$ has order $2,4$ or $6$, $D=11C$ works, since $\F_{11}$ does not contain any primitive third or fourth root of unity. 

When $K$ is real, we may assume (after multiplying $z$ with a rational number) that $z=(z_0,z_f)$ with $z_f \in \hat{\Z}^{\times}$ and $z_0 \in \R^{\times}_+$, and thus that its image in $(\OO_K \otimes \hat{\Z})^{\times}$ is congruent to a totally positive unit modulo $D$. It is thus enough to find $D$ such that any totally positive unit $u \in \OO_K^{\times}$ congruent to an integer modulo $D$ is congruent to $1$ modulo $C$.  

Fix a generator $u$ of the group of totally positive units of $\OO_K$ and suppose that there is no $D$ satisfying the previous condition. Then, with $D_n=(n+C)!$, there is an integer $d_n$ such that $u^{d_n}$ is congruent to an integer modulo $D_n$ but not to $1$ modulo $C$. Since $u$ is a totally positive unit, its conjugate under $\mrm{Gal}(K/\Q)$ is $u^{-1}$: thus $u^{2d_n}$ is congruent to $1$ modulo $D_n$. Let $\delta \in \hat{\Z}$ be any limit point of the $d_n$, then $u^{2\delta}=1$ in $(\OO_K \otimes \hat{\Z})^{\times}$. By \cite[Lemma 8.2.4]{CNBF}, one has $2\delta=0$ hence $\delta=0$. This implies that, for some $n$, $d_n$ is divisible by $|(\OO_K/C\OO_K)^{\times}|$, so that $u^{d_n}$ is congruent to $1$ modulo $C$, a contradiction.  

In either case, for a suitable value of $D$, the extension $K'$ is such that $\ker{\chi}$ contains $\ker{V_{K'}}$. This implies that there is a character $\psi_{0}: \mrm{Gal}(K'/K) \rar \C^{\times}$ such that $\psi_0 \circ V_{K'} = \chi$. In particular, if we view $\psi_{0}$ as a character $G_K^{ab} \rar \mrm{Gal}(K'/K) \rar \C^{\times}$, one has $\psi_0 \circ V=\chi$. \medskip

\emph{Step 2: Construction of infinitely many representations.} 

Let $c$ denote the action of the non-trivial element in $\mrm{Gal}(K/\Q)$ on $G_K^{ab}$ and fix an odd prime $p$ unramified in $K$. 

Let $S$ be any finite set of finite places of $K$. By Chebotarev, there exists an odd prime $q \equiv -1\pmod{p}$ inert in $K$ such that $q \notin S$. By local class field theory \cite[Theorems V.1.3, V.1.4]{NeuANT}, there exists a surjective character $\theta_q: G_{K_q} \rar \Z/p\Z$ whose restriction to the inertia subgroup has order $p$; in particular, if $\theta'_q$ denotes the conjugate of $\theta_q$ by the action of $\mrm{Gal}(K_q/\Q_q)$, $\theta'_q$ and $-\theta_q$ have the same restriction to the inertia group, so that the restriction to the inertia group of $\theta'_q-\theta_q$ has order $p$. By the Grunwald--Wang theorem \cite[Theorem 9.2.3]{CNBF}, there exists a continuous character $\theta: G_K \rar \Z/p\Z$ extending $\theta_q$: in particular, the restriction to the inertia group at $q$ of $\lambda := e^{\frac{2i\pi}{p}(\theta^c-\theta)}$ has order $p$. Moreover, one has $\lambda^p=\lambda\lambda^c=1$.  

By repeating the above construction, one can construct an increasing sequence $(q_i)_{i \geq 1}$ of odd primes and a sequence $(\lambda_i)_{i \geq 1}$ of continuous complex characters of $G_K$ such that:
\begin{itemize}[noitemsep,label=\tiny$\bullet$]
\item for every $i \geq 1$, $q_i \equiv -1\pmod{p}$ is inert in $K$ and unramified for $\psi_0,\lambda_1,\ldots,\lambda_{i-1}$, 
\item for every $i \geq 1$, $\lambda_i$ is ramified at $q_i$ and one has $\lambda_i^p=\lambda_i\lambda_i^c=1$, 
\end{itemize}

Define, for every $i \geq 1$, the character $\alpha_i=\psi_0\lambda_i$. Because $\lambda_i\lambda_i^c=1$ for each $i$, one has $\alpha_i\circ V = \psi_0\circ V=\chi$, so the Artin representation $\rho_i := \mrm{Ind}_K^{\Q}{\alpha_i}$ has determinant $\varepsilon \chi$, so it is in particular odd. Note that $\rho_i$ is unramified at $q_j$ for every $j > i$. 

Let $i \geq 1$ and $I_{q_i}$ be an inertia subgroup of $G_K$ at $q_i$. Then $\lambda_i(I_{q_i})$ is a subgroup of $\lambda_i(G_K) \simeq \Z/p\Z$ which is not reduced to the unit element, so $\lambda_i(I_{q_i})=\lambda_i(G_K)$ is the set of $p$-th roots of unity in $\C$. Therefore, the character $\alpha_i^{-1}\alpha_i^c=\lambda_i^{-2}\psi_0^{-1}\psi_0^c$ is ramified at $q_i$. In particular, $\alpha_i^{-1}\alpha_i^c$ is not trivial, so, by Mackey's criterion \cite[Corollary to Proposition 23]{Serre-Lin}, $\rho_i$ is irreducible.

Furthermore, since the restriction of $\alpha_i^{-1}\alpha_i^c$ to $I_{q_i}$ is not trivial, $(\rho_i)_{|I_{q_i}}$ is not given by scalar matrices. Hence, any twist of $\rho_i$ by a Dirichlet character is ramified at $q_i$, so, for any $1 \leq j < i$, $\rho_j$ is not the twist of $\rho_i$ by a Dirichlet character, and we are done.

}

\section{Examples and counter-examples}
\label{sect-counterex}
So far, we have identified in Theorem \ref{weakerconds} and in Section \ref{sect-theoA} some sufficient conditions to ensure a positive answer to Question \ref{loefflerqn}. On the other hand, the situation studied in Section \ref{sect-theoD} provides in some cases families of counter-examples. The purpose of this section is to study some situations where these results do not apply. We provide a few counter-examples which do not stem from Section \ref{sect-theoD}. We also showcase additional examples proving that, while the exceptional cases highlighted in Propositions \ref{special2} and \ref{specialq} can occur, they are not necessarily incompatible with Euler-adapted Galois images. Some of these counter-examples are explicit: as in the examples of Section \ref{sect-theoD}, we will then always refer to newforms by their LMFDB label \cite{lmfdb}. While this label only refers to the equivalence class of the modular forms under the action of $G_{\Q}$, the arguments work regardless of the choice of embedding of the modular forms in $\C$.

\prop[special-counter]{ Assume that $\mathfrak{p}$ is a good prime, and that $(f,g,p)$ falls in one of the following cases:
\begin{enumerate}[label=(\alph*),noitemsep]
\item\label{A5det} $f$ is $961.2.a.b$, $g$ is $3875.1.d.a$, and $p\equiv 13,37,83,107 \pmod{120}$.
\item\label{S4det} $f$ is $1849.2.a.g$, $g$ is $688.1.b.b$, and $p \equiv 7,17 \pmod{24}$. 
\item\label{A4det} $f$ is $1849.2.a.g$, $g$ is $2107.1.b.b$, and $p \equiv 11,13 \pmod{24}$. 
\item\label{S3det-CM} $f$ is $63.2.a.b$, $g$ is $1452.1.e.d$ and $p \equiv 5,7 \pmod{12}$.
\item\label{S3exc} $f$ is $189.2.a.f$, $g$ is $3468.1.i.a$ and $p \equiv \pm 5, \pm 43, \pm 67 \pmod{168}$. 
\end{enumerate}
Then $T_{f,g,\mathfrak{p}}$ does not satisfy (wE).}

\rem{In all these counter-examples, $f$ is a non-CM newform of weight $2$ with trivial character with coefficients in a quadratic number field, $g$ is always a newform of weight $1$ with quadratic character, except in case \ref{S3exc} where its character has order $4$. In particular, $\varepsilon_f\varepsilon_g \neq 1$. Moreover, $g$ never has complex multiplication by $\varepsilon_g\varepsilon_f$, so these counter-examples are not covered by Lemma \ref{efeg=1} or the construction of Section \ref{sect-theoD}. Finally, the field extension $H/\Q$ (associated to $f$ in the notation of Section \ref{sect-setup}) is quadratic. The corresponding situations are:
\begin{enumerate}[label=(\alph*),noitemsep]
\item $\tilde{\rho}_g(G_{\Q}) \simeq A_5$, $\varepsilon_g(G_H)=1$.
\item $\tilde{\rho}_g(G_{\Q}) \simeq S_4$, $\varepsilon_g(G_H)=1$.
\item $\tilde{\rho}_g(G_{\Q}) \simeq A_4$, $\varepsilon_g(G_H)=1$.
\item $\varepsilon_g(G_H)=1$ and $g$ has CM by a field distinct from $H$.
\item $i \in \varepsilon_g(G_H)$, $g$ has CM by a field distinct from $H$, $\tilde{\rho}_g(G_{\Q}) = \tilde{\rho}_g(G_H)$ is a dihedral group of order $6$. 
\end{enumerate}}

\demo{In this proof, $j$ (in any domain of characteristic zero) denotes a choice of primitive third root of unity. Moreover, if $R$ is a ring and $\alpha \in R^{\times}$, $\GL{R}_{\det=\alpha}$ denotes the subset of matrices in $\GL{R}$ with determinant $\alpha$. 

The proof strategy is the same in all cases: using the LMFDB data, we check that in every case, $F=\Q$, $L,H$ are quadratic number fields, $f$ has trivial character, and $\mfk{p}$ is good and coprime to $|\rho_g(G_{\Q})|$. Then, for each coset $C=\sigma G_{H(\mu_{p^{\infty}})} \subset G_{\Q(\mu_{p^{\infty}})}$, we can determine an element $\alpha \in \OO_{L,\mfk{p}}^{\times}$ attached to $C$ by Corollary \ref{papier}, since we know that $\alpha^2 \in F$. 

Assume that there is $\tau \in \sigma G_{\Q(\mu_{p^{\infty}})}$ such that the kernel of the action of $\tau-1$ on $T_{f,g,\mfk{p}}$ has rank one. By Corollary \ref{papier}, we can find matrices $N \in \rho_g(\sigma)\rho_g(G_H)$ and $M \in \begin{pmatrix}\alpha & 0\\0 & \alpha^{-1}\end{pmatrix}\SL{\Z_p}$ such that the kernel of $M \otimes N -1$ is a line. This condition implies that $\det{N}=\det{M}\det{N} \neq 1$ by Lemma \ref{efeg=1}, that neither $M$ and $N$ is scalar, and that $N$ has an eigenvalue $u$ such that $u^{-1}$ is an eigenvalue of $M$. Note that $u$ is a root of unity with order $\omega$ prime to $p$.   

Now, $\alpha^2 \in \Z_p^{\times}$, so 
\[\begin{pmatrix}\alpha & 0\\0 & \alpha^{-1}\end{pmatrix}\SL{\Z_p} = \alpha^{-1}\begin{pmatrix}\alpha^2 & 0\\0 & 1\end{pmatrix}\SL{\Z_p} = \alpha^{-1}\GL{\Z_p}_{\det=\alpha^2}.\]

Then $\alpha u^{-1}$ is an eigenvalue of the matrix $M_1 = \alpha M \in \GL{\Z_p}_{\det=\alpha^2}$. Therefore \[\mrm{Tr}(M_1) = \alpha u^{-1}+\frac{\alpha^2}{\alpha u^{-1}} = \alpha(u^{-1}+u) \in \Z_p.\]

Thus either $\alpha u, \alpha u^{-1}$ both lie in $\Z_p^{\times}$, or they are Galois conjugates of each other lying in some quadratic extension of $\Q_p$. Since $p$ is odd and prime to the order of $u$, and since $\alpha^2 \in \Z_p^{\times}$, $\alpha u, \alpha u^{-1}$ are contained in the maximal unramified extension of $\Q_p$, so the absolute arithmetic Frobenius $\varphi$ of $\Q_p^{unr}$ either fixes them both or exchanges them. When $\alpha \notin \Q_p$, since $\alpha^2 \in \Q$, one has $\varphi(\alpha)=-\alpha$, hence $\varphi(\alpha u)=-\alpha u^p$. We thus reach a contradiction if $-1 \notin \{u^{p-1},u^{p+1}\}$. When $\alpha \in \Q_p$, one has $\varphi(\alpha u)=\alpha u^p$, so we reach a contradiction if $1 \notin \{u^{p-1},u^{p+1}\}$. 

Assuming that $\varepsilon_f$ is the trivial character, that $F=\Q$, and that $\mfk{p}$ is coprime to $2|\rho_g(G_{\Q})|$, we have thus proved the following. Pick, for every coset $C$ of $G_{H(\mu_{p^{\infty}})}$ in $G_{\Q(\mu_{p^{\infty}})}$, some $\alpha_C \in \OO_{L,\mfk{p}}^{\times}$ given by Corollary \ref{papier}, and let $V_C \subset L_0$ be a set containing the eigenvalues of non-scalar matrices $N \in \rho_g(C)$ such that $\det{N} \neq 1$. If, for each $C$ such that $\alpha_C \notin \Z_p^{\times}$, $V_C$ does not contain an element $u$ such that $-1 \in \{u^{p-1},u^{p+1}\}$, and, for each $C$ such that $\alpha_C \in \Z_p^{\times}$, $V_C$ does not contain an element $u$ such that $1 \in \{u^{p-1},u^{p+1}\}$, then $T_{f,g,\mfk{p}}$ does not satisfy (wE). 

This is the property that we will use in all five cases.\medskip

\ref{A5det} According to the LMFDB, $\varepsilon_f$ is trivial, $\varepsilon_g$ is the unique quadratic character of conductor $31$, $L=\Q(\sqrt{2})$, $f$ is not CM and has an inner twist by $\varepsilon_g$, so that $F=\Q$ and $H=\Q(\sqrt{-31})$. Since $\varepsilon_g(G_H)=1$, Lemma \ref{efeg=1} shows that any $\sigma \in G_{\Q(\mu_{p^{\infty}})}$ such that the kernel of $\sigma-1$ has rank one cannot lie in $G_{H}$. 

Let $D$ be the derived subgroup of $\rho_g(G_{\Q})$ and $Z \leq \rho_g(G_{\Q})$ be the subgroup of scalar matrices. Since $g$ is of $A_5$-type, $D \subset \SL{L_0}$ is a finite group whose projective image is $A_5$. In particular, it has even order, hence contains $-I_2$. Thus we have an exact sequence $1 \rar \{\pm I_2\} \rar D \rar A_5 \rar 1$. Moreover, $D$ and $\rho_g(G_{\Q})$ have the same projective image, so $\rho_g(G_{\Q})=ZD$. Since $\det{\rho_g(G_{\Q})}=\{\pm 1\}$ and $D \subset \SL{L_0}$, one has $Z=\{\pm iI_2, \pm I_2\}$, so that $\rho_g(G_{\Q})=D \cup iD$. Since $D \subset \rho_g(G_H) \subset \rho_g(G_{\Q})$ and $\varepsilon_g(G_H)=1$, one has $\rho_g(G_{H})=D$, and $\rho_g(G_{\Q}\backslash G_H)=iD$.

Because $\varepsilon_g(G_H)=\{1\}$, it is enough to test the coset $G_{\Q(\mu_{p^{\infty}})}\backslash G_{H(\mu_{p^{\infty}})}$. The element $\alpha \in L^{\times}$ attached to this coset by Corollary \ref{papier} is $\alpha=\sqrt{2}$. Because $p \equiv 3,5 \pmod{8}$, $\alpha \notin \Q_p$. Hence, we need to prove that if $u$ is an eigenvalue of any non-scalar matrix in $iD$, one does not have $-1 \in \{u^{p-1},u^{p+1}\}$. Let $N \in iD$ be a non-scalar matrix. Then $\det{N}=-1$, and its projective image $\tilde{N}$ has order $2$, $3$ or $5$. 

\begin{itemize}[noitemsep,label=$-$]
\item if $\tilde{N}$ has order $2$, then $N^2$ is scalar, so Cayley--Hamilton implies that $N^2=I_2$, so the eigenvalues of $N$ are $\pm 1$, so $u^{p \pm 1}=1 \neq -1$.  
\item if $\tilde{N}$ has order $3$, then $N^3$ is scalar. Since $\det{N}=-1$, the eigenvalues of $N$ are $\{ij,ij^2\}$ or $\{-ij,-ij^2\}$. Either way, any eigenvalue $u$ of $N$ has order $12$, and $p \equiv \pm 1 \pmod{12}$, so $-1 \notin \{u^{p-1},u^{p+1}\}$. 
\item if $\tilde{N}$ has order $5$, then $N^5$ is scalar. Since $\det{N}=-1$, any eigenvalue $u$ of $N$ is a root of unity of order $20$. Since $p \not\equiv \pm 1 \pmod{5}$, the order of $u^{p \pm 1}$ is divisible by $5$, so $-1 \notin \{u^{p-1},u^{p+1}\}$ and we are done. 
\end{itemize}

\ref{S4det} According to the LMFDB, in this case, $\varepsilon_f$ is trivial, $L=\Q(\sqrt{6})$, the projective image of $\rho_g$ is $S_4$ and $\varepsilon_g$ is the unique non-trivial quadratic character of conductor $43$. Moreover, $f$ has no CM and has an inner twist by $\varepsilon_g$, so that $F=\Q$ and $H=\Q(\sqrt{-43})$. 

Let $D=\rho_g(G_{\Q})'$ and $Z$ denote the group of scalar matrices contained in $\rho_g(G_{\Q})$. One has $\det{Z} \subset \varepsilon_g(G_{\Q}) = \{\pm 1\}$ so $Z \subset \{\pm I_2, \pm i I_2\}$. Moreover, the coefficient field of $g$ is $\Q(\sqrt{-2})$, so that $g$ has no coefficient equal to $\pm 2i$, hence $Z\subset \{\pm I_2\}$. As in the previous case, $D \subset \rho_g(G_{\Q}) \cap \SL{L_0}$ and its projective image is $S_4'=A_4$, so $D$ has even order, hence $-I_2 \in D$ and $Z=Z \cap D=\{\pm I_2\}$. Thus $D$ is a subgroup of $\rho_g(G_{\Q})$ with index $2$; since $\varepsilon_g(G_H)=\{1\}$, one furthermore has $D \subset \rho_g(G_H) \subset \rho_g(G_{\Q}) \cap \SL{L_0}$: thus $D=\rho_g(G_H)=\rho_g(G_{\Q}) \cap \SL{L_0}$ and $\rho_g(G_{\Q} \backslash G_H)=\rho_g(G_{\Q}) \cap \GL{L_0}_{\det=-1}$. 

As above, $\varepsilon_g(G_H)=\{1\}$, so it is enough to test the coset $C := G_{\Q(\mu_{p^{\infty}})}\backslash G_{H(\mu_{p^{\infty}})}$. The element $\alpha \in \OO_{L,\mfk{p}}^{\times}$ attached to $C$ by Corollary \ref{papier} is $\sqrt{6}$. By the previous paragraph, the matrices of $\rho_g(C)$ have determinant $-1$ and their projective images lie in $S_4 \backslash A_4$, so their orders are $2$ or $4$. Hence, for any non-scalar $M \in \rho_g(C)$, the set of eigenvalues of $M$ is either $\{1,-1\}$ or $\{\omega,-\omega^{-1}\}$ for some primitive $8$-th root of unity $\omega$. %

Because $p \equiv \pm 1 \pmod{8}$ and $p \not\equiv \pm 1 \pmod{12}$, $2$ is a square in $\Z_p^{\times}$ but $3$ is not, so $\alpha \notin \Z_p^{\times}$. Thus, we need to prove that $-1$ is not contained in $\{1^{p\pm 1},(-1)^{p \pm 1},\omega^{p \pm 1}\}=\{1,\omega^{p \pm 1}\}$ where $\omega$ is a primitive $8$-th root of unity. Since $p \equiv \pm 1 \pmod{8}$, $p \pm 1$ is congruent modulo $8$ to $0,2,-2$, so $\{\omega^{p \pm 1}\} \subset \{1,\omega^2,\omega^{-2}\}=\{1,\pm i\}$ does not contain $-1$, and we are done.

\medskip

\ref{A4det} According to the LMFDB, $L,F,H$ and $\varepsilon_f,\varepsilon_g$ are as above (so $\varepsilon_g$ is the quadratic character of $H/\Q$), and $\tilde{\rho}_g(G_{\Q}) \simeq A_4$, while $\rho_g(G_{\Q})$ has cardinality $48$: hence the scalar matrices contained in $\rho_g(G_{\Q})$ are exactly $\{\pm iI_2, \pm I_2\}$. Now, $\rho_g(G_H)$ is a subgroup of index at most $2$ of $\rho_g(G_{\Q})$; since it is contained in $\SL{L_0}$, it does not contain $iI_2$, so $\rho_g(G_{\Q} \backslash G_H) = i\rho_g(G_H)$. In particular, $\tilde{\rho}_g(G_H)=\tilde{\rho}_g(G_{\Q})$.%

Since $\varepsilon_g(G_H)=\{1\}$, it is enough to consider the coset $G_{\Q(\mu_{p^{\infty}})}\backslash G_{H(\mu_{p^{\infty}})}$, to which the element $\alpha$ attached by Corollary \ref{papier} is $\sqrt{6}$. Since $p \not\equiv \pm 1 \pmod{8}$ and $p \equiv \pm 1 \pmod{12}$, $2$ is not a square in $\Z_p^{\times}$ and $3$ is, so $\alpha \notin \Q_p$. 

Hence we need to show that for any non-scalar matrix $M \in \rho_g(G_{\Q} \backslash G_H)=i\rho_g(G_H)$, for any eigenvalue $u$ of $M$, one has $-1 \notin \{u^{p-1},u^{p+1}\}$. %

Let $N \in i\rho_g(G_H)$ be a non-scalar matrix, then $\det{N}=-1$ and its image $\tilde{N} \in \PGL{L_0}$ has order $2$ or $3$. 
\begin{itemize}[noitemsep,label=$-$]
\item If $\tilde{N}$ has order $2$, then $N \sim \Delta_{1,-1}$, in which case the claim is clear. %
\item If $\tilde{N}$ has order $3$, then as in case \ref{A5det}, any eigenvalue $u$ of $N$ has order $12$. Since $p \equiv \pm 1\pmod{12}$, one among $u^{p-1}$ and $u^{p+1}$ is equal to $1$, and the other has order $6$, so neither equals $-1$ and we are done. %
\end{itemize}

\ref{S3det-CM} In this case, $L=\Q(\sqrt{3})$, $F=\Q$, $H=\Q(j)$, $\varepsilon_f$ is trivial and $\varepsilon_g$ is the character of $H/\Q$, whereas $g$ has complex multiplication by $\Q(\sqrt{-11})$. Because $\varepsilon_g$ is the character of the imaginary quadratic field $H/\Q$, one has $G_H=\rho_g^{-1}(\SL{L_0})$. Moreover, the field of coefficients of $g$ is $\Q(j)$, so the group $Z$ of scalar matrices contained in $\rho_g(G_{\Q})$ has order dividing $6$: since $\varepsilon_g(G_{\Q})=\{\pm 1\}$, $Z \subset \{\pm I_2\}$. Therefore, the morphism $\det: \tilde{\rho}_g(G_{\Q}) \rar \{\pm 1\}$ is well-defined and surjective with kernel $\tilde{\rho}_g(G_H)$. %

Because $\varepsilon_g(G_H)=\{1\}$, we only need to consider the coset $C := G_{\Q(\mu_{p^{\infty}})} \backslash G_{H(\mu_{p^{\infty}})}$. An $\alpha \in \OO_{L,\mfk{p}}^{\times}$ attached to $C$ by Corollary \ref{papier} is $\alpha=\sqrt{3}$. Since $p \equiv 5,7 \pmod{12}$, $3$ is not a square in $\Q_p$, so we need to show that if $u$ is any eigenvalue of some non-scalar matrix in $\rho_g(C)=\rho_g(G_{\Q}) \backslash \SL{L_0}=\rho_g(G_{\Q})\backslash \rho_g(G_H)$, then $-1 \notin \{u^{p-1},u^{p+1}\}$.   

The group $\tilde{\rho}_g(G_H)$ is a subgroup of index $2$ of the dihedral group $D := \tilde{\rho}_g(G_{\Q})$ of order $12$, and it is not contained in the maximal cyclic subgroup of $D$, so the elements of $\rho_g(G_{\Q}) \backslash \rho_g(G_H)$ fall in two categories: matrices $N$ with $\det{N}=-1$ and projective order $2$, and matrices $N$ with $\det{N}=-1$ and projective order $6$.

\begin{itemize}[noitemsep,label=$-$]
\item If $\det{N}=-1$ and the projective image of $N$ has order $2$, its eigenvalues are $\pm 1$; one has $1^{\pm 1}=(-1)^{p \pm 1}=1 \neq -1$, so we are done. 
\item If $\det{N}=-1$ and the projective image of $N$ has order $6$, and $u$ is an eigenvalue of $N$, then $N \sim \Delta_{u,-u^{-1}}$, $u^6=u^{-6}$, $u^3 \neq -u^{-3}$, $u^2 \neq u^{-2}$. Hence $u^2 \neq 1$, $u^6 = \pm 1$ and $u^6 \neq -1$, so $u^6=1$. Therefore, $u^{3(p \pm 1)}=1$, so $-1 \notin \{u^{p \pm 1}\}$.
\end{itemize} 

\ref{S3exc} In this case, $L=\Q(\sqrt{7}),F=\Q$ and $H=\Q(\sqrt{-3})$; $g$ has complex multiplication by $\Q(\sqrt{-51})$, $\varepsilon_g$ is a character of conductor $51$ and order $4$, and $\tilde{\rho}_g(G_{\Q}) \simeq S_3$; $\varepsilon_f$ is the trivial character. Since $\varepsilon_g(2893)=\pm i$, $\varepsilon_g(G_H)=\{\pm 1, \pm i\}=\varepsilon_g(G_{\Q})$. The LMFDB shows that $\rho_g(G_{\Q})$ is isomorphic to $S_3 \times \Z/8\Z$. Its center (which is scalar, by Schur's lemma) is thus the set of $8$-th roots of unity, and $\rho_g(G_{\Q}) \cap \SL{L_0} \rar \tilde{\rho}_g(G_{\Q})$ is onto with kernel $\{\pm I_2\}$. 

By the LMFDB data, one has $G_H \subsetneq G_H\cdot \ker{\tilde{\rho}_g}$, so $\tilde{\rho}_g(G_H) = \tilde{\rho}_g(G_{\Q})$. Using MAGMA \cite{magma}, we can also see that $\ker{\rho_g} \subset G_H$, so that $[\rho_g(G_{\Q}):\rho_g(G_H)]=2$. 

For the sake of easier notation, write $C \subset \rho_g(G_{\Q})$ for the pre-image of the cyclic subgroup of order $3$ of $\tilde{\rho}_g(G_{\Q})$. 

Since $\rho_g(G_H)$ and $\rho_g(G_{\Q})$ have the same determinant and projective image, we can split elements of $\rho_g(G_{\Q})$ in six categories: $\{\pm i I_2, \pm I_2\}$, scalar matrices of order $8$, non-scalar matrices in $C$ with determinant $\pm 1$, non-scalar matrices in $C$ with determinant $\pm i$, matrices outside $C$ with determinant $\pm 1$, matrices outside $C$ with determinant $\pm i$. Then $\rho_g(G_H)$ is exactly made with the matrices of the first, third, and last categories. 

An element $\alpha \in \OO_{L,\mfk{p}}^{\times}$ attached to any $\sigma \in G_{\Q(\mu_{p^{\infty}})} \backslash G_H$ (resp. any $\sigma \in G_{H(\mu_{p^{\infty}})}$) by Corollary \ref{papier} is $\alpha=\sqrt{7}$ (resp. $\alpha=1$). Note that because $p \equiv \pm 5, \pm 43, \pm 67 \pmod{168}$, $\sqrt{7} \notin \Q_p$. 

Therefore, we need to show that:
\begin{itemize}[noitemsep,label=$-$]
\item If $u$ is any eigenvalue of any matrix $N$ in the third or last category such that $\det{N} \neq 1$, then $p$ is not congruent to $\pm 1$ modulo the order of $u$.
\item If $u$ is any eigenvalue of any matrix $N$ in the fourth or fifth category such that $\det{N} \neq 1$, then $-1 \notin \{u^{p-1},u^{p+1}\}$. 
\end{itemize}  

Let $N \in C$ be in the third category with $\det{N} \neq 1$. Then $\det{N} = -1$ and $N^3$ is scalar with determinant $-1$, so as in \ref{A5det} any eigenvalue $u$ of $N$ has order $12$. We are done since $p \not\equiv \pm 1 \pmod{12}$. %

Let $N$ be a matrix in the sixth category. Then the projective image of $N$ has order $2$ and $\det{N}=\pm i$, so $N^2=\pm iI_2$, and any eigenvalue $u$ of $N$ is a primitive eighth root of unity. We are done since $p \not\equiv \pm 1 \pmod{8}$. %

Let $N \in C$ be a matrix in the fourth category. Then $\det{N}=\pm i$ and $N^3$ is scalar. Then $N \sim \omega \Delta_{j,j^2}$ for some primitive $8$-th root of unity $\omega$, so any eigenvalue $u$ of $N$ has order $24$. Since $p \not\equiv \pm 1 \pmod{12}$, one has $1 \notin \{u^{2(p\pm 1)}\}$, so $-1 \notin \{u^{p \pm 1}\}$, and we are done. %

Let $N$ be a matrix in the fifth category such that $\det{N} \neq 1$. Then $\det{N}=-1$ and the projective image of $N$ has order $2$, and we conclude as in \ref{A5det}. 

}

The proposition below, on the other hand, shows that, even in the case where the assumption of Proposition \ref{special2} is satisfied (and therefore the reasoning of Corollary \ref{sufficient-for-euler-adapted} does not apply), the representation can still have Euler-adapted Galois image.

\prop[pvaries-fine]{Let $f \in \mathcal{S}_k(\Gamma_1(N)), g \in \mathcal{S}_1(\Gamma_1(M))$ be newforms with $k \geq 2$ and respective characters $\varepsilon_f, \varepsilon_g$. Assume that $f$ has no complex multiplication and let $H, L, F$ be the fields attached to $f$ as in Section \ref{sect-setup}. Let $L_0 \subset \C$ be a number field containing $L$, the Fourier coefficients of $g$ and the eigenvalues of all matrices in $\rho_g(G_{\Q})$. Assume furthermore that:
\begin{itemize}[noitemsep,label=\tiny$\bullet$]
\item $g$ has CM by an imaginary quadratic field $K \not\subset H$,
\item $\rho_g(G_{\Q})$ contains a scalar matrix with order $4$,  
\item $\tilde{\rho}_g(G_{\Q})$ has cardinality $4n$ for some odd integer $n > 1$,
\item $[\rho(G_{\Q}):\rho(G_H)]=[\tilde{\rho}_g(G_{\Q}):\tilde{\rho}_g(G_H)]=2$, 
\item $\varepsilon_g(G_{\Q})=\varepsilon_g(G_H)=\{\pm 1, \pm i\}$, 
\item $\rho_g(\ker{\varepsilon_f}) \neq \rho_g(G_H)$.
\end{itemize}
Then, for all but finitely many good primes $\mfk{p} \subset L_0$ for $(f,g)$, the representation $T_{f,g,\mfk{p}}$ satisfies (sE). 
}

\demo{If $\rho_g(G_H)$ contains a matrix with eigenvalues $1$ and $-1$, we are done by Corollary \ref{sufficient-for-euler-adapted}. We may thus assume that $\rho_g(G_H)$ does not contain any matrix whose eigenvalues are $1$ and $-1$. 

Suppose that $\varepsilon_g(G_{HK})=\{\pm 1,\pm i\}$. Then $\varepsilon_g(G_H \backslash G_K)$ is also equal to $\{\pm 1, \pm i\}$, so we can find $\sigma \in G_H \backslash G_K$ such that $\varepsilon_g(\sigma)=-1$. Because $g$ has CM by $K$, $\rho_g(\sigma)$ has trace zero, so by Cayley--Hamilton $\rho_g(\sigma)$ is similar to $\Delta_{1,-1}$, which contradicts our assumption. 

Since $G_H/G_{HK}$ is a group of order $2$, it follows that $\varepsilon_g(G_{HK})=\{\pm 1\}$, so $\varepsilon_g(G_H \backslash G_K)=\{\pm i\}$. Let $H_0 \subset H$ be the subfield such that $\ker{\varepsilon_f}=G_{H_0}$. One has $\rho_g(G_H) \subset \rho_g(G_{H_0}) \subset \rho_g(G_{\Q})$, the first inclusion is proper and $[\rho_g(G_{\Q}):\rho_g(G_H)]=2$ by assumption, so $H \neq H_0$ and $\rho_g(G_{\Q})=\rho_g(G_{H_0})$. In particular, $\rho_g(G_{H_0})$ contains the image $N_-$ of the complex conjugation, which is not contained in $\rho_g(G_H)$ by assumption, since its eigenvalues are $1,-1$.  

The quotient map $\rho_g(G_{\Q})/\rho_g(G_H) \rar \tilde{\rho}_g(G_{\Q})/\tilde{\rho}_g(G_H)$ is surjective between groups of same cardinality, so it is injective. Therefore, $\rho_g(G_{\Q})$ and $\rho_g(G_H)$ have the same scalar matrices. Any pre-image by $\rho_g$ of any scalar matrix lies in $G_K$, thus $\rho_g(G_{\Q})$ and $\rho_g(G_{HK})$ contain the same scalar matrices. Because $i \notin \varepsilon_g(G_{HK})$, $\rho_g(G_{HK})$ does not contain a scalar matrix with order $8$, so neither does $\rho_g(G_{\Q})$. Thus, the scalar matrices in $\rho_g(G_{\Q})$ are exactly $\pm I_2, \pm iI_2$. 

Next, we show that there is $N_4 \in \rho_g(G_{H_0}) \backslash \rho_g(G_H) = \rho_g(G_{\Q}) \backslash \rho_g(G_H)$ whose image in $\PGL{L_0}$ has order $2$ such that $\det{N_4}$ has order $4$. Indeed, let $\delta: \tilde{\rho}_g(G_{\Q}) \rar \varepsilon_g(G_{\Q})/\{\pm 1\}$ be the determinant map: it is onto a group of cardinality $2$. Then $\delta$ vanishes on the complex conjugation, so $\delta$ does not vanish on the cyclic subgroup $\tilde{\rho}_g(G_K)$ of $\tilde{\rho}_g(G_{\Q})$. This cyclic subgroup has order $2n$ (and $n$ is odd), so it contains a unique element $\kappa$ of order $2$, and $\delta(\kappa)$ is non-trivial. Because $K \not\subset H$ and $\tilde{\rho}_g(G_H) \neq \tilde{\rho}_g(G_{\Q})$, $\tilde{\rho}_g(G_K) \cap \tilde{\rho}_g(G_H)$ is the subgroup of index $2$ of $\tilde{\rho}_g(G_K)$, so it has odd cardinality and does not contain $\kappa$. Thus, if $N_4 \in \rho_g(G_{\Q})$ is any element lifting $\kappa$, then $N_4 \notin \rho_g(G_H)$, the projective image of $N_4$ has order $2$, and $\det{N_4}$ has order $4$.  

By the second part of Proposition \ref{jointimage-H}, we can fix an element $\tau \in G_{H_0(\mu_{p^{\infty}})}$ such that $\rho_g(\tau) \notin \rho_g(G_H)$. Let $\mfk{p} \subset \OO_{L_0}$ be a good prime for $(f,g)$ of residue characteristic $p$. Let $\alpha \in  \OO_{L,\mfk{p}}^{\times}$ be the element attached to $\tau$ by Corollary \ref{papier}. Then $\alpha^2 \in \OO_{F,\mfk{p}}^{\times}$ and $\alpha \notin F$. We will show that the representation $T_{f,g,\mfk{p}}$ satisfies (sE) by distinguishing between the following three cases:

\begin{enumerate}[label=(\roman*),noitemsep]
\item\label{pvaries-1} if $2$ is a square in $F_{\mfk{p}}$,
\item\label{pvaries-2} if $\alpha^2 \in F$ is a square in $F_{\mfk{p}}$,
\item\label{pvaries-3} if $2\alpha^2$ is a square in $F_{\mfk{p}}$.
\end{enumerate}

This is enough to conclude because $2,\alpha^2$ are units in $F_{\mfk{p}}$.  

In case \ref{pvaries-1}, pick $\sigma_1 \in G_H \backslash G_K$. Since $g$ has CM by $K$, $\rho_g(\sigma_1)$ has trace zero and determinant $\pm i$, so its eigenvalues are $u,-u$ for some primitive eighth root of unity $u \in L_0$. We can pick a square root $t$ of $2$ in $L_0$ such that $tu-1$ is a fourth root of unity, then $t \in (L_0)_{\mfk{p}}$ is an element of its subfield $F_{\mfk{p}}$. By Proposition \ref{jointimage-H}, we can choose $\sigma \in G_{H(\mu_{p^{\infty}})}$ such that $\rho_g(\sigma)=\rho_g(\sigma_1)$ and the characteristic polynomial of $\sigma$ acting on $T_{f,\mfk{p}}$ is given by $x^2-tx+1$. Then the eigenvalues of $\sigma$ acting on $T_{f,g,\mfk{p}}$ (with ring of coefficients $\OO_{(L_0)_{\mfk{p}}}$) are given by the quadruple $(u\cdot u, u^{-1}u, u(-u), u^{-1}(-u))=(u^2,1,-u^2,-1)$, so $1$ is a simple eigenvalue for the actions of the action of $\sigma$ on $T_{f,g,\mfk{p}}$ and $T_{f,g,\mfk{p}}/\mfk{p}$, so it fulfills the condition (sE).

In case \ref{pvaries-2}, by Corollary \ref{papier}, we may pick $\sigma \in \tau G_{H(\mu_{p^{\infty}})}$ such that $\sigma$ acts on $T_{f,\mfk{p}}$ with matrix $U=\begin{pmatrix}1 & 1\\0 & 1\end{pmatrix}$ and $\rho_g(\sigma)=N_-$, and the conclusion follows. 

In case \ref{pvaries-3}, let $u \in L_0$ be a primitive eighth root of unity such that $u,-u$ are the two eigenvalues of $N_4$. Then the two elements $\alpha u+\alpha u^{-1} = \pm \sqrt{2} \alpha, (\alpha u)(\alpha u^{-1})=\alpha^2$ of $(L_0)_{\mfk{p}}$ lie in the subring $\OO_{F_{\mfk{p}}}$, so there exists a matrix $M \in \GL{\OO_{F_{\mfk{p}}}}$ with eigenvalues $\alpha u, \alpha u^{-1}$, so that $\det{M}=\alpha^2$. By Corollary \ref{papier}, there exists $\sigma \in \tau G_{H(\mu_{p^{\infty}})}$ acting on $T_{f,\mfk{p}}$ by the matrix $\alpha^{-1}M$ (in a suitable basis) and such that $\rho_g(\sigma)=N_4$. Then the eigenvalues of $\sigma$ acting on $T_{f,g,\mfk{p}}$ are given as in case \ref{pvaries-1} by $u^2, 1, -u^2, -1$, all of which are distinct mod $\mfk{p}$, and we are done.    

}

\cor[special2-exc-nothold]{Let $f \in \mathcal{S}_2(\Gamma_0(63)), g \in \mathcal{S}_1(\Gamma_1(675))$ be newforms with respective labels $63.2.a.b$ and $675.1.g.a$. Then $T_{f,g,\mathfrak{p}}$ has Euler-adapted Galois image for every good prime $\mathfrak{p}$. }

\demo{By Proposition \ref{only-we-se}, we are reduced to studying the condition (sE).

In this case, the LMFDB says that $L=\Q(\sqrt{3}), F=\Q, H=\Q(\sqrt{-3})$, $f$ does not have complex multiplication, and $\varepsilon_g$ is a character of conductor $5$ and order $4$. Moreover, $\tilde{\rho}_g(G_{\Q})$ is the dihedral group of order $12=4 \cdot 3$, and $g$ has complex multiplication by $\Q(\sqrt{-15})$, which is not contained in $H$. Furthermore, the number field LMFDB data shows that $\ker{\tilde{\rho}_g} \subset G_H$, so $[\tilde{\rho}_g(G_{\Q}):\tilde{\rho}_g(G_H)]=[\rho_g(G_{\Q}):\rho_g(G_H)]=2$. Since $a_{139}(g)=\pm 2i$, $iI_2 \in \rho_g(G_H) \subset \rho_g(G_{\Q})$. Since the coefficient field of $g$ is $\Q(i,\sqrt{6})$, it does not contain a primitive $8$-th root of unity, so the scalar matrices of $\rho_g(G_{\Q})$ are exactly $\pm iI_2, \pm I_2$. Finally, $\varepsilon_g(G_H)=\varepsilon_g(G_{\Q})=\{\pm i, \pm 1\}$, since $\varepsilon_g$ and $H$ have coprime conductors. 

Therefore we can apply Proposition \ref{pvaries-fine}.}

\rem{In the situation of Corollary \ref{special2-exc-nothold}, the assumptions of Proposition \ref{special2} hold, and those of Proposition \ref{specialq} do not. Hence none of the sufficient conditions of Corollary \ref{sufficient-for-euler-adapted} applies. 

In fact, the following also holds: there is no couple \[(M,N) \in \left(\SL{\Q} \times \rho_g(G_H)\right) \cup \left(\begin{pmatrix}\sqrt{3} & 0\\0 & 1/\sqrt{3}\end{pmatrix}\SL{\Q} \times (\rho_g(G_{\Q}) \backslash \rho_g(G_H))\right)\] such that the kernel of $M \otimes N -1$ is a line. However, what the proof of Corollary \ref{special2-exc-nothold} shows is that this claim holds after replacing the coefficient field $\Q$ with its completion at any place $v \neq 2,3$. 
}

~\\

\footnotesize{ 

\nocite{*}
\bibliographystyle{amsplain}
\bibliography{prod-final}
}

\medskip
\medskip
\noindent 

\medskip
\medskip
\noindent
\footnotesize{\textsc{Universit\'e Paris Cit\'e and Sorbonne Universit\'e, CNRS, IMJ-PRG, F-75013 Paris, France.}\\\emph{Email adress: }{\sffamily e.n.g.studnia@math.leidenuniv.nl}}

\end{document}